\documentclass[final,leqno]{siamltex}
\usepackage{amsfonts}
\usepackage{amssymb}
\usepackage{amsmath}
\usepackage{amscd}
\usepackage{url}
\usepackage{graphicx}
\usepackage{psfrag,epsfig}
\usepackage{exscale}

\newcommand{\dx}{\,\mathrm{d}x}
\newcommand{\dX}{\,\mathrm{d}X}
\newcommand{\px}{\partial x}
\newcommand{\pX}{\partial X}
\newcommand{\R}{\mathbb{R}}
\newcommand{\M}{\mathbb{M}}
\newcommand{\Rr}{R}
\newcommand{\dxi}[1]{\partial_{x_{#1}}}

\newcommand{\Eq}{}

\newcommand{\Nedelec}{N\'ed\'elec}

\newcommand{\contraPiola}{\mathcal{F}^{\mathrm{div}}}
\newcommand{\coPiola}{\mathcal{F}^{\mathrm{curl}}}

\newcommand{\triang}{\mathcal{T}}

\DeclareMathOperator{\Div}{div}
\DeclareMathOperator{\Curl}{curl}
\DeclareMathOperator{\Rot}{rot}
\DeclareMathOperator{\Grad}{grad}
\DeclareMathOperator{\ExtD}{d}
\DeclareMathOperator{\AdjExtD}{\delta}

\DeclareMathOperator{\skw}{skw}
\DeclareMathOperator{\tr}{tr}

\newcommand{\flux}{\sigma}
\newcommand{\stress}{\sigma}

\newcommand{\HdivO}{H(\mathrm{div}; \Omega)}
\newcommand{\HcurlO}{H(\mathrm{curl}; \Omega)}
\newcommand{\Hdiv}{H(\mathrm{div})}
\newcommand{\Hcurl}{H(\mathrm{curl})}

\newcommand{\RT}{\mathrm{RT}}
\newcommand{\BDM}{\mathrm{BDM}}
\newcommand{\BDFM}{\mathrm{BDFM}}
\newcommand{\NED}{\mathrm{NED}}
\newcommand{\DG}{\mathrm{DG}}
\newcommand{\Poly}{\mathcal{P}}

\newtheorem{thm}{Theorem}[section]

\newtheorem{remark}[thm]{{\it Remark}}
\newtheorem{example}[thm]{{\it Example}}
\newtheorem{defi}[thm]{Definition}

\usepackage{verbatim}
\newenvironment{code}[1]%
{\center\tabular{c}\hline\\ \footnotesize\minipage{#1\textwidth}\verbatim}
{\endverbatim\endminipage\\ \\ \hline\endtabular\endcenter}

\title{Efficient Assembly of {\boldmath$H(\mathrm{\lowercase{div}})$} and {\boldmath$H(\mathrm{\lowercase{curl}})$} Conforming Finite Elements\thanks{Received by the editors October 23, 2008;
accepted for publication (in revised form) September 9, 2009; published electronically November 20, 2009. \URL sisc/31-6/73901.html}}

\author{Marie E. Rognes\thanks{Centre of Mathematics for Applications, University of Oslo, P.O.~Box~1053, 0316 Oslo, Norway (meg@cma.uio.no).}
  \and Robert C. Kirby\thanks{Department of Mathematics and Statistics, Texas Tech University, P.O.~Box~1042, Lubbock, TX\@,
  79409-1042 (robert.c.kirby@ttu.edu). This author's work was supported by the United States Department of Energy Office of Science under grant DE--FG02--07ER25821.}
  \and Anders Logg\thanks{Center for Biomedical Computing, Simula Research Laboratory,
    Department of Informatics, University of Oslo, P.O.~Box~134, 1325 Lysaker, Norway (logg@simula.no).
    This author's work was supported by a Center of Excellence grant from the Research Council of Norway
    to the Center for Biomedical Computing at Simula Research
    Laboratory and by
    an Outstanding Young Investigator grant from the Research Council of Norway, NFR 180450.}}

\begin{document}
\slugger{sisc}{2009}{31}{6}{4130--4151}
\maketitle

\setcounter{page}{4130}

\begin{abstract}
  In this paper, we discuss how to efficiently evaluate and assemble
  general finite element variational forms on $H(\mathrm{div})$ and $H(\mathrm{curl})$.
  The proposed strategy relies on a decomposition of the element
  tensor into a precomputable reference tensor and a mesh-dependent
  geometry tensor. Two key points must then be considered: the
  appropriate mapping of basis functions from a reference element, and
  the orientation of geometrical entities. To address these issues, we
  extend here a previously presented representation theorem for
  affinely mapped elements to Piola-mapped elements. We also discuss a
  simple numbering strategy that removes the need to contend with
  directions of facet normals and tangents.  The result is an
  automated, efficient, and easy-to-use implementation that allows a
  user to specify finite element variational forms on $H(\mathrm{div})$ and
  $H(\mathrm{curl})$ in close to mathematical notation.
\end{abstract}



\begin{DOI}
10.1137/08073901X
\end{DOI}

\pagestyle{myheadings}
\thispagestyle{plain}
\markboth{MARIE E. ROGNES, ROBERT C. KIRBY, AND ANDERS LOGG}%
     {ASSEMBLY OF $\Hdiv$ AND $\Hcurl$ FINITE ELEMENTS}
\section{Introduction}

The Sobolev spaces $\Hdiv$ and $\Hcurl$ play an important role in many
applications of mixed finite element methods to partial differential
equations. Examples include second order elliptic partial differential
equations, Maxwell's equations for electromagnetism, and the linear
elasticity equations. Mixed finite element methods may provide
advantages over standard $H^1$ finite element discretizations in terms
of added robustness, stability, and flexibility. However, implementing
$\Hdiv$ and $\Hcurl$ methods requires additional code complexity for
constructing basis functions and evaluating variational forms, which
helps to explain their relative scarcity in practice.

The FEniCS project~\cite{logg:www:03, Log2007a} comprises a
collection of free software components for the automated solution
of differential equations. One of these components is the FEniCS
form compiler (FFC)~\cite{KirLog2006, KirLog2007a, logg:www:04}.
FFC allows finite element spaces over simplicial meshes and
multilinear forms to be specified in a form language close to the
mathematical abstraction and notation. The form compiler generates
low-level (C++) code for efficient form evaluation and assembly
based on an efficient tensor contraction. Moreover, the FErari
project~\cite{KirKne2005, Kirby:2008:BDC, KirLogEtAl2006,
KirSco07} has developed specialized techniques for further
optimizing this code based on underlying discrete structure. FFC
relies on the FInite element Automatic Tabulator
(FIAT)~\cite{Kir04, www:FIAT, Kir06} for the tabulation of finite
element basis functions. FIAT provides methods for efficient
tabulation of finite\break element basis functions and their
derivatives at any particular point. In particular, FIAT provides
simplicial $\Hdiv$ element spaces such as the families of
Raviart and Thomas~\cite{RavTho77b},
Brezzi, Douglas, and Marini~\cite{BreDou85}, and
Brezzi et al.\ \cite{BreDou87}, as well as
$\Hcurl$ elements of the \Nedelec{} types~\cite{Ned80, Ned86}.

\enlargethispage*{7pt}
Previous iterations of FFC have enabled easy use of $H^1$ and $L^2$
conforming finite element spaces, including discontinuous Galerkin
formulations, but support for $\Hdiv$ and $\Hcurl$ spaces has been
absent. In this paper, we extend the previous work~\cite{KirLog2006,
KirLog2007a, OelLog2008a} to allow simple and efficient compilation of
variational forms on $\Hdiv$ and $\Hcurl$, including mixed
formulations on combinations of $H^1$, $\Hdiv$, $\Hcurl$, and
$L^2$. The efficiency of the proposed approach relies, in part, on the
tensor representation framework established in~\cite{KirLog2007a}. In
this framework, the element tensor is represented as the contraction
of a reference tensor and a geometry tensor. The former can be
efficiently precomputed given automated tabulation of finite element
basis functions on a reference element, while the latter depends on
the geometry of each physical element. For this strategy, a key aspect
of the assembly of $\Hdiv$ and $\Hcurl$ conforming element spaces
becomes the Piola transformations, isomorphically mapping basis
functions from a reference element to each physical element. Also, the
orientation of geometrical entities such as facet tangents and normals
must be carefully considered.

Implementations of $\Hdiv$ and $\Hcurl$ finite element spaces, in
particular of arbitrary degree, are not prevalent. There are, to our
knowledge, no implementations that utilize the compiled approach to
combine the efficiency of low-level optimized code with a fully
automated high-level interface. Some finite element packages, such as
FEAP~\cite{www:FEAP}, do not provide $\Hdiv$ or $\Hcurl$ type elements
at all. Others, such as FreeFEM~\cite{www:FreeFEM}, typically provide
only low-order elements such as the lowest-order Raviart--Thomas
elements. Some libraries such as deal.II~\cite{www:deal.II} or
FEMSTER~\cite{CasRie05} do provide arbitrary degree elements of
Raviart--Thomas and \Nedelec{} type, but do not automate the evaluation
of variational forms. NGSolve~\cite{www:ngsolve} provides arbitrary
order $\Hdiv$ and $\Hcurl$ elements along with automated assembly, but
only for a predefined set of bilinear forms.

This exposition and the FFC implementation consider the assembly of
$\Hdiv$ and $\Hcurl$ finite element spaces on simplicial
meshes. However, the underlying strategy is extendible to
nonsimplicial meshes and tensor-product finite element spaces defined
on such meshes. A starting point for an extension to isoparametric $H^1$
conforming finite elements was discussed in~\cite{KirLog2006}. The
further extensions to $\Hdiv$ and $\Hcurl$ follow the same lines as
for the simplicial case discussed in this note.

The outline of this paper is as follows. We begin by reviewing basic
aspects of the function spaces $\Hdiv$ and $\Hcurl$ in
section~\ref{sec:hdivcurl}, and we provide examples of variational forms
defined on these spaces.  We continue, in section~\ref{sec:elements},
by summarizing the $\Hdiv$ and $\Hcurl$ conforming finite elements
implemented by FIAT. In section~\ref{sec:multilinearforms}, we recap
the multilinear form framework of FFC and we present an extension of the
representation theorem from~\cite{KirLog2007a}. Subsequently, in
section~\ref{sec:assembling}, we provide some notes on the assembly of
$\Hdiv$ and $\Hcurl$ elements. Particular emphasis is placed on
aspects not easily found in the standard literature, such as the choice of
orientation of geometric entities.  In section~\ref{sec:examples}, we
return to the examples introduced in section~\ref{sec:hdivcurl} and
illustrate the ease and terseness with which even complicated mixed
finite element formulations may be expressed in the FFC form
language. Convergence rates in agreement with theoretically predicted
results are presented to substantiate the veracity of the
implementation. Finally, we make some concluding remarks in
section~\ref{sec:conclusion}.

\section{\boldmath$\Hdiv$ and $\Hcurl$}
\label{sec:hdivcurl}

In this section, we summarize some basic facts about the Sobolev
spaces $\Hdiv$ and $\Hcurl$ and we discuss conforming finite element
spaces associated with them. Our primary focus is on properties
relating to interelement continuity and change of variables. The
reader is referred to the monographs~\cite{BreFor91} and~\cite{Mon03}
for a more thorough analysis of $\Hdiv$ and $\Hcurl$, respectively.

\subsection{Definitions}

For an open domain $\Omega \subset \R^n$, we let $L^2(\Omega, \R^n)$
denote the space of square-integrable vector fields on $\Omega$ with
the associated norm $||\cdot||_{0}$ and inner-product $\langle \cdot,
\cdot \rangle$, and we abbreviate $L^2(\Omega) = L^2(\Omega, \R^1)$. We
define the following standard differential operators on smooth fields
$v$: $D^{\alpha} v = \dxi{1}^{\alpha_1} \cdots \dxi{m}^{\alpha_{m}} v$
for a multi-index $\alpha$ of length $m$, $ \Div v = \sum_{i=1}^n
\dxi{i} v_i$, $\Curl v = (\dxi{2} v_3 - \dxi{3} v_2, \dxi{3} v_1 -
\dxi{1} v_3, \dxi{1} v_2 - \dxi{2} v_1)$, and $\Rot v = \dxi{1} v_2 -
\dxi{2} v_1$.  We may then define the spaces $H^m(\Omega)$, $\HdivO$,
and $\HcurlO$ by
\begin{align*}
  H^m(\Omega) &= \{ v \in L^2(\Omega): \; D^{\alpha} v \in
  L^2(\Omega), \; |\alpha| \leq m\}, \quad m = 1, 2, \dots, \\ \HdivO
  &= \{ v \in L^2(\Omega, \R^n): \; \Div v \in L^2(\Omega) \}, \\
  \HcurlO &= \left \{
  \begin{tabular}{ll}
    $\{ v \in L^2(\Omega, \R^2): \; \Rot v  \in L^2(\Omega) \}$,
    & $\Omega \subset \R^2$, \\
    $\{ v \in L^2(\Omega, \R^3): \; \Curl v  \in L^2(\Omega, \R^3) \}$,
    & $\Omega \subset \R^3$,
  \end{tabular}
  \right .
\end{align*}
with derivatives taken in the distributional sense. The reference to
the domain $\Omega$ will be omitted when appropriate, and the
associated norms will be denoted $||\cdot||_{m}$, $||\cdot||_{\Div}$,
and $||\cdot||_{\Curl}$.  Furthermore, we let $\M$ denote the space of
matrices and we let $H(\Div; \Omega, \M)$ denote the space of
square-integrable matrix fields with square-integrable row-wise
divergence.

{\spaceskip .26em plus .1em minus .1em For the sake of compact notation, we shall also adopt the exterior
calculus notation of~\cite{ArnFalWin06} and let $\Lambda^k(\Omega)$
denote the space of smooth differential $k$-forms on $\Omega$, and let
$L^2 \Lambda^k(\Omega)$ denote the space of square-integrable
differential $k$-forms on $\Omega$.  We further let $\ExtD$ denote the
exterior derivative with adjoint $\AdjExtD$, and we define $H
\Lambda^k(\Omega) = \{ v \in L^2 \Lambda^k(\Omega), \ExtD v \in L^2
\Lambda^k(\Omega) \}$. Further, $\mathcal{P}_r \Lambda^k$ is the space
of polynomial $k$-forms of up to and including degree $r$, and
$\mathcal{P}^{-}_r \Lambda^k$ denotes the reduced space as defined in
\cite[section~3.3]{ArnFalWin06}.}

\subsection{Examples}

The function spaces $\Hdiv$ and $\Hcurl$ are the natural function
spaces for an extensive range of partial differential equations, in
particular in mixed formulations. We sketch some examples in the
following, both for motivational purposes and for later reference. The
examples considered here are mixed formulations of the Hodge Laplace
equations, the standard eigenvalue problem for Maxwell's equations, and
a mixed formulation for linear elasticity with weakly imposed
symmetry. We return to these examples in section \ref{sec:examples}.

\begin{example}[mixed formulation of Poisson's equation]
\label{ex:mixedpoisson}{\rm
The most immediate example involving the space $\Hdiv$ is a mixed
formulation of Poisson's equation: $ -\Delta u = f$ in $\Omega \subset
\R^n$. By introducing the flux $\flux = -\Grad u$ and assuming
Dirichlet boundary conditions for $u$, we obtain the following mixed
variational problem: Find $\sigma \in \HdivO$ and $u \in L^2(\Omega)$
satisfying
\begin{equation}
  \label{eq:mixedpoisson}
  \langle \tau, \flux \rangle - \langle \Div \tau, u \rangle
  + \langle v, \Div \flux  \rangle = \langle v, f \rangle
\end{equation}
for all $\tau \in \HdivO$ and $v \in L^2(\Omega)$.}
\end{example}
\begin{example}[the Hodge Laplacian]
\label{ex:hodgelaplace}
{\rm With more generality, we may consider weak formulations of the
Hodge Laplacian equation $(\ExtD \AdjExtD + \AdjExtD \ExtD) u = f$ on a
domain $\Omega \subset \R^n$; see \cite[section 7]{ArnFalWin06}. For
simplicity of presentation, we assume that $\Omega$ is contractible
such that the space of harmonic forms on $\Omega$ vanishes. The
formulation in Example~\ref{ex:mixedpoisson} is the equivalent of
seeking $u \in H \Lambda^n$ and $\flux = \AdjExtD u \in H
\Lambda^{n-1}$ for $n = 2, 3$ with natural boundary conditions (the
appropriate trace being zero).  To see this, we test $\flux = \AdjExtD
u$ against $\tau \in H \Lambda^{n-1}$ and we test $(\ExtD \AdjExtD +
\AdjExtD \ExtD) u = f$ against $v \in H \Lambda^n$ to obtain}
\begin{displaymath}
  \langle \tau, \flux \rangle - \langle \tau , \AdjExtD u \rangle +
  \langle v, \ExtD \sigma \rangle =
  \langle v, f \rangle,
\end{displaymath}
{\rm noting that $d u = 0$ for $u \in H \Lambda^n$.
Integrating by parts, we obtain}
\begin{equation} \label{eq:mixedpoisson,d}
  \langle \tau, \flux \rangle - \langle \ExtD \tau , u \rangle +
  \langle v, \ExtD \sigma \rangle =
  \langle v, f \rangle.
\end{equation}
{\rm We may restate~(\ref{eq:mixedpoisson,d}) in the
form~(\ref{eq:mixedpoisson}) by making the identifications
$\AdjExtD u = - \Grad u$, $\ExtD \tau = \Div \tau$, and $\ExtD \sigma
 = \Div \sigma$.
If $\Omega \subset \R^3$, we may also consider the following mixed
formulations of the Hodge Laplace equation.}
\renewcommand{\theenumi}{\roman{enumi}}
\begin{enumerate}
\item[{\rm (i)}] {\rm Find} $\flux \in H \Lambda^1 = \Hcurl$ {\rm and} $u \in H \Lambda^2 =
\Hdiv$ {\rm such that}
\begin{equation}
  \label{eq:curldiv}
  \langle \tau, \flux  \rangle
  - \langle \Curl \tau, u \rangle
  + \langle v, \Curl \flux \rangle
  + \langle \Div v, \Div u \rangle = \langle v, f \rangle
\end{equation}
{\rm for all} $\tau \in H \Lambda^1$, $v \in H \Lambda^2$.
\item[{\rm (ii)}] {\rm Find} $\flux \in H \Lambda^0 = H^1$ {\rm and} $u \in H \Lambda^1 =
\Hcurl$ {\rm such that}
\begin{equation}
  \label{eq:gradcurl}
  \langle \tau, \flux \rangle
  - \langle \Grad \tau, u \rangle
  + \langle v, \Grad \flux \rangle
  + \langle \Curl v, \Curl u \rangle = \langle v, f \rangle
\end{equation}
{\rm for all $\tau \in H \Lambda^0, v \in H \Lambda^1$.}
\end{enumerate}
\end{example}

\begin{example}[cavity resonator]
  \label{ex:cavity}{\rm
  The time-harmonic Maxwell equations in a cavity with perfectly
  conducting boundary induces the following eigenvalue problem: Find
  resonances $\omega \in \R$ and eigenfunctions $E \in H_0(\Curl;
  \Omega)$,  satisfying}
  \begin{equation}
    \label{eq:maxwell}
    \langle \Curl F, \Curl E \rangle = \omega^2 \langle F, E \rangle
    \quad \forall F \in H_0(\Curl; \Omega),
  \end{equation}
  {\rm where $H_0(\Curl; \Omega) = \{v \in \HcurlO | \; v \times n
  |_{\partial \Omega}= 0 \}$. Note that the formulation (\ref{eq:maxwell})
  disregards the original divergence-free constraint for the electric
  field $E$ and thus includes the entire kernel of the $\Curl$
  operator, corresponding to $\omega = 0$ and electric fields of the
  form $E = \Grad \psi$.}
\end{example}
\begin{example}[elasticity with weakly imposed symmetry]
  \label{ex:elasticity}{\rm
  Navier's equations for linear elasticity can be reformulated using
  the stress tensor $\stress$, the displacement $u$, and an additional
  Lagrange multiplier $\gamma$ corresponding to the symmetry of the
  stress constraint. The weak equations for $\Omega \subset \R^2$,
  with the natural\footnote{Note that the natural boundary condition
  in this mixed formulation is a Dirichlet condition, whereas for
  standard $H^1$ formulations the natural boundary condition would be
  a Neumann\break condition.} boundary condition $u|_{\partial \Omega} = 0$,
  take the following form: Given $f \in L^2(\Omega, \R^n)$, find
  $\stress \in H(\Div; \Omega, \M)$, $u \in L^2(\Omega,
  \R^n)$, and $\gamma \in L^2(\Omega)$ such that
  \begin{equation}
    \label{eq:elasticity}
    \langle \tau, A \stress \rangle + \langle \Div \tau, u \rangle +
    \langle v, \Div \stress \rangle + \langle \skw \tau, \gamma \rangle +
    \langle \eta, \skw \stress \rangle = \langle v, f \rangle
  \end{equation}
  for all $\tau \in H(\Div; \Omega, \M)$, $v \in
  L^2(\Omega, \R^n)$, and $\eta \in L^2(\Omega)$. Here, $A$ is the
  compliance tensor, and $\skw \tau$ is the scalar representation of
  the skew-symmetric component of $\tau$; more precisely, $2 \skw
  \tau = \tau_{21} - \tau_{12}$. This formulation has the advantage of
  being robust with regard to nearly incompressible materials and it provides an alternative foundation for complex materials with
  nonlocal stress-strain relations. For more details, we refer the
  reader to \cite{ArnFalWin07}.}
\end{example}

\subsection{Continuity-preserving mappings for \boldmath$\Hdiv$ and $\Hcurl$}
At this point, we turn our attention to a few results on
continuity-preserving mappings for $\Hdiv$ and $\Hcurl$. The results
are classical and we refer the reader to~\cite{BreFor91, Mon03} for a more
thorough treatment.

First, it follows from Stokes' theorem that in order for piecewise
$\Hdiv$ vector fields to be in $\Hdiv$ globally, the traces of the
normal components over patch interfaces must be continuous, and
analogously tangential continuity is required for piecewise $\Hcurl$
fields. More precisely, we have the following: Let $\triang_h = \{K\}$
be a partition of $\Omega$ into subdomains. Define the space
$\Sigma_h$ of piecewise $\Hdiv$ functions relative to this partition
$\triang_h$:
\begin{equation}
  \Sigma_h = \{ \phi \in L^2(\Omega, \R^n): \; \phi|_K \in
  H(\Div; K)\ \forall K \in \triang_h\}.
\end{equation}
Then $\phi \in \Sigma_h$ is in $\HdivO$ if and only if the normal traces
of $\phi$ are continuous across all element interfaces. Analogously, if
$\phi|_{K} \in H(\Curl; K)$ for all $K \in \triang_h$,
then $\phi \in \HcurlO$ if and only if the tangential traces are
continuous across all element interfaces.

Second, we turn to consider a nondegenerate mapping $F: \Omega_0
\rightarrow F(\Omega_0) = \Omega$ with Jacobian $DF(X)$, $X \in
\Omega_0 \subset \R^n$. For $\Phi \in H^m(\Omega_0)$, the mapping
$\mathcal{F}$ defined by
\begin{equation}
  \label{eq:affine}
  \mathcal{F}(\Phi) = \Phi \circ F^{-1}
\end{equation}
is an isomorphism from $H^m(\Omega_0)$ to $H^m(\Omega)$. This,
however, is not the case for $\Hdiv$ or $\Hcurl$, since $\mathcal{F}$
does not in general preserve continuity of normal or tangential
traces. Instead, one must consider the contravariant and covariant
Piola mappings which preserve normal and tangential continuity,
respectively.

\begin{defi}[the contravariant and covariant Piola mappings]
  \label{def:Piolas}
  Let $\Omega_0 \subset \R^n$, let $F$ be a nondegenerate mapping
  from $\Omega_0$ onto $F(\Omega_0) = \Omega$ with $J = DF(X)$, and let
  $\Phi \in L^2(\Omega_0, \R^n)$.

  The contravariant Piola mapping $\contraPiola$ is defined by
  \begin{equation}
    \label{eq:contraPiola}
    \contraPiola(\Phi) = \frac{1}{\det J} J \Phi \circ F^{-1}.
  \end{equation}
  The covariant Piola mapping $\coPiola$ is defined by
  \begin{equation}
    \label{eq:coPiola}
    \coPiola(\Phi) = J^{-T} \Phi \circ F^{-1}.
  \end{equation}
\end{defi}
\unskip%

\begin{remark}{\rm
We remark that the contravariant Piola mapping is usually defined with
an absolute value, $\contraPiola(\Phi) = \frac{1}{|\det J|} J \Phi
\circ F^{-1}$. However, omitting the absolute value, as
in~\eqref{eq:contraPiola}, can simplify the assembly of $\Hdiv$
elements, as will be expounded in section~\ref{sec:assembling}.}
\end{remark}

The contravariant Piola mapping is an isomorphism of $H(\Div;
\Omega_0)$ onto $\HdivO$, and the covariant Piola mapping is an
isomorphism of $H(\Curl; \Omega_0)$ onto $\HcurlO$. In particular, the
contravariant Piola mapping preserves normal traces, and the covariant
Piola mapping preserves tangential traces. We illustrate this below in
the case of simplicial meshes in two and three space dimensions
(triangles and tetrahedra). The same results hold for nonsimplicial
meshes with cell-varying Jacobians, such as quadrilateral
meshes~\cite{art:ArnBofFal2005, BreFor91}.

\begin{example}[Piola mapping on triangles in $\R^2$]{\rm
Let $K_0$ be a triangle with vertices $X^i$ and edges $E^i$ for
$i=1,2,3$. We define the unit tangents by $T^i = E^i / ||E^i||$. We
further define the unit normals by $N^i = \Rr T^i$, where
\begin{equation}
  \Rr =
  \begin{pmatrix}
    0 & 1 \\
    -1 & 0
  \end{pmatrix}
\end{equation}
is the clockwise rotation matrix.

Now, assume that $K_0$ is affinely mapped to a (nondegenerate)
simplex $K$ with vertices $x^i$. The affine mapping $F_K: K_0
\rightarrow K$ takes the form $x = F_K(X) = J X + b$ and satisfies
$x^i = F_K(X^i)$ for $i = 1, 2, 3$.  It follows that edges are mapped
by
\begin{align*}
  e = x^i - x^j = J (X^i - X^j) = J E.
\end{align*}
Similarly, normals are mapped by
\begin{equation*}
  ||e|| n
  = \Rr e
  = \Rr J E
  = (\det J) J^{-T} \Rr E
  = (\det J) J^{-T} ||E|| N,
\end{equation*}
where we have used that $\frac{1}{\det J} \Rr J \Rr^T = J^{-T}$
and thus $\Rr J = (\det J) J^{-T} \Rr$ for $J \in \mathbb{R}^{2 \times 2}$.

The relation between the mappings of tangents and normals (or edges and
rotated edges) may be summarized in the following commuting diagrams:
\begin{align}
  \label{cd:normalmap}
  \begin{CD}
    T @> J ||E|| / ||e|| >> t \\
    @VRVV @VVRV \\
    N @> (\det J) J^{-T} ||E|| / ||e|| >> n
  \end{CD}
  \quad\quad\quad
  \begin{CD}
    E @> J >> e \\
    @VRVV @VVRV \\
    ||E|| N @> (\det J) J^{-T} >> ||e|| n
  \end{CD}
\end{align}

\begin{figure}
  \begin{center}
{\hspace*{-4pc}{%
    \psfrag{0}{$1$}
    \psfrag{1}{$2$}
    \psfrag{2}{$3$}
    \psfrag{K0}{$K_0$}
    \psfrag{K}{$K$}
    \psfrag{Fn}{$\Phi_n(X) = (\frac{1}{\sqrt{2}},\frac{1}{\sqrt{2}})$}
    \psfrag{Ft}{$\Phi_t(X) = (-\frac{1}{\sqrt{2}},\frac{1}{\sqrt{2}})$}
    \psfrag{fn1}{$\phi_n(x) = (1,\frac{1}{2})$}
    \psfrag{ft1}{$\phi_t(x) = (0,\frac{1}{2})$}
    \psfrag{fn2}{$\phi_n(x) = (\frac{1}{2},0)$}
    \psfrag{ft2}{$\phi_t(x) = (-\frac{1}{2},1)$}
    \psfrag{x}{$x$}
    \psfrag{X}{$X$}
    \psfrag{Fdiv}{$\contraPiola$}
    \psfrag{Fcurl}{$\coPiola$}
    \includegraphics[width=0.88\textwidth]{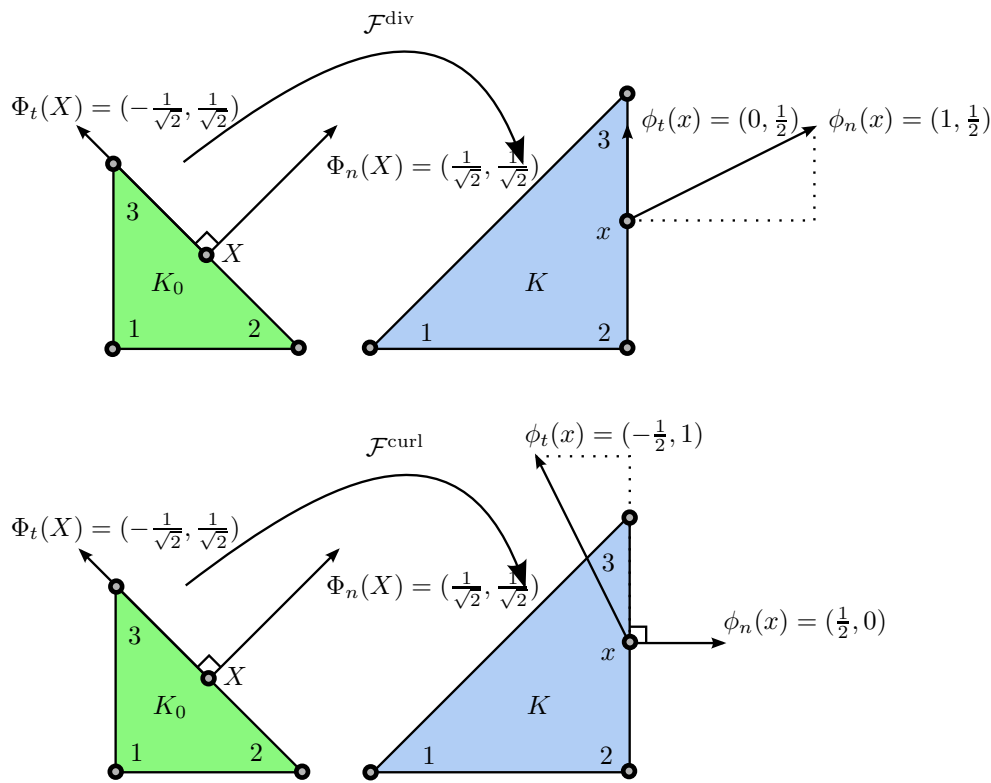}}}
    \caption{Mapping two vector fields $\Phi_n$ and $\Phi_t$ between
      two triangles using the contravariant and covariant Piola
      mappings.  The contravariant Piola mapping (above) preserves
      normal traces of vector fields, and the covariant Piola mapping
      (below) preserves tangential traces of vector fields. This means
      in particular that the contravariant Piola mapping maps tangents
      to tangents (which have a zero normal component), and that the
      covariant Piola mapping maps normals to normals (which have a
      zero tangential component). Note that this is somewhat
      counterintuitive; the contravariant $\Hdiv$ Piola mapping
     always maps tangential fields to tangential fields but does not
      in general map normal fields to normal fields. However, in both
      cases the normal component (being zero and one, respectively) is
      preserved.}
    \label{fig:piolaexample}
  \end{center}
\end{figure}

With this in mind, we may study the effect of the Piola transforms on
normal and tangential traces. Let $\Phi \in C^{\infty}(K_0, \R^n)$ and
let $\phi = \contraPiola(\Phi)$. Then
\begin{align*}
  ||e|| \, \phi(x) \cdot n
  =
  ||e||
  \left( (\det J)^{-1} J \Phi(X) \right)^T
  \left( (\det J) J^{-T} ||E||/||e|| N \right)
  =
  ||E|| \, \Phi(X) \cdot N.
\end{align*}
Thus, the contravariant Piola mapping preserves normal traces for
vector fields under affine mappings, up to edge lengths. In general,
the same result holds for smooth, nondegenerate mappings $F_K$ if the
Jacobian $DF_K(X)$ is invertible for all $X \in K_0$.

Similarly, let $\phi = \coPiola(\Phi)$. Then
\begin{equation}
  ||e|| \, \phi(x) \cdot t
  =
  ||e||
  \left( J^{-T} \Phi(X) \right)^T
  \left( J ||E|| / ||e|| \right)
  =
  ||E|| \, \Phi(X) \cdot T.
\end{equation}
Thus, the covariant Piola preserves tangential traces for vector
fields, again up to edge lengths. Observe that the same result holds
for tetrahedra without any modifications.
The effect of the contravariant and covariant Piola mappings on normal
and tangential traces is illustrated in Figure~\ref{fig:piolaexample},
where $||E|| = ||e||$ for simplicity.}
\end{example}

\begin{example}[contravariant Piola mapping on tetrahedra in $\R^3$]{\rm
  Now, let $K_0$ be a tetrahedron. As explained above, the covariant
  Piola mapping preserves tangential traces. To study the effect of
  the contravariant Piola mapping on normal traces, we define the face
  normals of $K$ by $N = \frac{E^i \times E^j}{||E^i \times E^j||}$.
  Then
  \begin{align*}
    ||e^i \times e^j|| n
    = J E^i \times J E^j
    = \det J J^{-T} (E^i \times E^j)
    = ||E^i \times E^j|| \det J J^{-T} N,
  \end{align*}
  since $(J u) \times (J v) = \det J J^{-T} (u \times v)$. Let $\Phi
  \in C^{\infty}(K_0, \R^n)$ and let $\phi =
  \contraPiola(\Phi)$. Then, it follows that
  \begin{equation*}
    ||e^i \times e^j|| \, \phi(x) \cdot n
    =
    ||E^i \times E^j|| \, \Phi(X) \cdot N.
  \end{equation*}
  Thus, the contravariant Piola mapping preserves normal traces, up to
  the area of faces.}
\end{example}

We finally remark that if $J \in \R^{2 \times 2}$ defines a conformal,
orientation-preserving map, the contravariant and covariant Piola
mappings coincide. In $\R^3$, $J$ must also be orthogonal for this to
occur.

\section{\boldmath$\Hdiv$ and $\Hcurl$ conforming finite elements}
\label{sec:elements}

To construct $\Hdiv$ and $\Hcurl$ conforming finite element spaces,
that is, discrete spaces $V_h$ satisfying $V_h \subset \Hdiv$ or $V_h
\subset \Hcurl$, one may patch together local function spaces (finite
elements) and make an appropriate matching of degrees of freedom over
shared element facets. Here, a facet denominates any geometric entity
of positive codimension in the mesh (such as an edge of a triangle or
an edge or face of a tetrahedron). In particular, one requires that
degrees of freedom corresponding to normal traces match for $\Hdiv$
conforming discretizations and that tangential traces match for
$\Hcurl$ conforming discretizations.

Several families of finite element spaces with degrees of freedom
chosen to facilitate this exist. For $\Hdiv$ on simplicial
tessellations in two dimensions, the classical conforming families
are those of Raviart and Thomas ($\RT_r$, $r = 0, 1, 2,
\dots$)~\cite{RavTho77b}; Brezzi, Douglas, and Marini ($\BDM_r$,
$r = 1, 2, \dots$)~\cite{BreDou85}; and Brezzi et al.\ ($\BDFM_r$, $r = 1, 2, \dots$)~\cite{BreDou87}. The
former two families were extended to three dimensions by
\Nedelec~\cite{Ned80, Ned86}. However, the same notation will be
used for the two- and three-dimensional $\Hdiv$ element spaces
here. For $\Hcurl$, there are the families of \Nedelec{} of the
first kind ($\NED_r^1$, $r = 0, 1, 2, \dots$)~\cite{Ned80} and of
the second kind ($\NED_r^2$, $r = 1, 2, \dots$)~\cite{Ned86}. We
summarize in Table~\ref{tab:listofelements} those $\Hdiv$ and
$\Hcurl$ conforming finite elements that are supported by FIAT and
hence by FFC. In general, FFC can wield any finite element space
that may be generated from a local basis through either of the
aforedescribed mappings. In Table~\ref{tab:approximability}, we
also summarize some basic approximation properties of these
elements for later comparison with numerical results in
section~\ref{sec:examples}.

\begin{table}[hbt]
\renewcommand{\arraystretch}{1.4}
\footnotesize
\caption{$\Hdiv$ and $\Hcurl$ conforming finite elements on
    triangles and tetrahedra supported by FIAT and FFC for $r \geq
    1$. When applicable, the elements are listed with their exterior
    calculus notation, along with their original references. Note that
    for $K \subset\R^3$, the Raviart--Thomas and Brezzi--Douglas--Marini
    elements are also known as the first and second kind $\Hdiv$
    \Nedelec{} (face)\break elements, respectively.}
  \label{tab:listofelements}
\centering
 \begin{tabular}{|c||c|c|}
      \hline
      Simplex & $\Hdiv$ & $\Hcurl$ \\
      \hline \hline
      $K \subset \R^2$
      &
      \begin{tabular}{p{1.5cm}p{1.5cm}p{0.6cm}}
        $\mathrm{BDM}_r$  & $\Poly_{r} \Lambda^1(K)$      & \cite{BreDou85}  \\
        $\mathrm{RT}_{r-1}$   & $\Poly_{r}^{-} \Lambda^1(K)$ & \cite{RavTho77b} \\
        $\mathrm{BDFM}_r$ & ---                           & \cite{BreDou87} \\
      \end{tabular}
      &
      \begin{tabular}{p{1.6cm}p{1.5cm}p{0.6cm}}
          $\mathrm{NED}^1_{r-1}$ & --- &
      \end{tabular}
      \\ \hline
      $K \subset \R^3$
      &
      \begin{tabular}{p{1.5cm}p{1.5cm}p{0.6cm}}
        $\mathrm{BDM}_r$ & $\Poly_{r} \Lambda^2(K)$      & \cite{Ned86} \\
        $\mathrm{RT}_{r-1}$  & $\Poly_{r}^{-} \Lambda^2 (K)$ & \cite{Ned80} \\
        $\mathrm{BDFM}_r$ & ---                         &  \\
      \end{tabular}
      &
      \begin{tabular}{p{1.5cm}p{1.5cm}p{0.6cm}}
        & & \\
        $\mathrm{NED}^1_{r-1}$ & $\Poly_{r}^- \Lambda^1 (K)$    & \cite{Ned80} \\
        \\
      \end{tabular}
      \\ \hline
    \end{tabular}
\vspace*{-6pt}
\end{table}

\begin{table}[hbt]
\renewcommand{\arraystretch}{1.4}
\footnotesize
\caption{Approximation properties of the spaces from Table
    {\rm \ref{tab:listofelements}}. $C > 0$, $r \geq 1$, $1 \leq m \leq
    r$. $\Pi_h$ denotes the canonical interpolation operator, defined
    by the degrees of freedom, onto the element space in question. For
    simplicity, it is assumed that $v$ is sufficiently smooth for the
    interpolation operators to be well-defined, and for the given
    norms to be bounded. For more details and sharper estimates; cf.~{\rm \cite{BreFor91, Mon03}}.}
  \label{tab:approximability}
  \centering
    \begin{tabular}{|l|ll|l|}
      \hline
      Finite element &
      Interpolation estimates &   \\
      \hline \hline
      $\Poly_r \Lambda^k(\Omega)$
      & $||v - \Pi_h v||_0 \leq C h^{m+1} ||v||_{m+1}$,
      & $||v - \Pi_h v||_{\Div, \Curl} \leq C h^{m} ||v||_{m+1}$ \\
      \hline
      $\Poly_r^{-} \Lambda^k(\Omega)$
      & $||v - \Pi_h v||_0 \leq C h^{m} ||v||_{m}$,
      &  $||v - \Pi_h v||_{\Div, \Curl} \leq C h^{m} ||v||_{m+1}$ \\
      \hline
      $\BDFM_r$
      & $||v - \Pi_h v||_0 \leq C h^{m} ||v||_{m}$,
      & $||v - \Pi_h v||_{\Div} \leq C h^{m} ||v||_{m+1}$ \\
      \hline
    \end{tabular}
\end{table}

For the reasons above, it is common to define the degrees of freedom
for each of the elements in Table~\ref{tab:listofelements} as moments
of either normal or tangential traces over element facets. However,
one may alternatively consider point values of traces at suitable
points on element facets (in addition to any internal degrees of
freedom). Thus, the degrees of freedom for the lowest order
Raviart--Thomas space on a triangle may be chosen as the normal
components at the edge midpoints, and for the lowest order
Brezzi--Douglas--Marini space, we may consider the normal components
at two points on each edge (positioned symmetrically on each edge and
not touching the vertices). This, along with the appropriate scaling
by edge length, is how the degrees of freedom are implemented in FIAT.

\section{Representation of \boldmath$\Hdiv$ and $\Hcurl$ variational forms}
\label{sec:multilinearforms}

In this section, we discuss how multilinear forms on $\Hdiv$ or
$\Hcurl$ may be represented as a particular tensor contraction,
allowing for precomputation of integrals on a reference element and
thus efficient assembly of linear systems. We follow the notation
from~\cite{KirLog2006, KirLog2007a} and extend the representation
theorem from~\cite{KirLog2007a} for multilinear forms on $H^1$ and
$L^2$ to $\Hdiv$ and $\Hcurl$. The main new component is that we must
use the appropriate Piola mapping to map basis functions from the
reference element.

\subsection{Multilinear forms and their representation}

Let $\Omega \subset \R^n$ and let $\{V_h^j\}_{j=1}^{\rho}$ be a set of
finite dimensional spaces associated with a tessellation $\triang =
\{K\}$ of $\Omega$. We consider the following canonical linear
variational problem: Find $u_h \in V_h^2$ such that
\begin{equation} \label{eq:varproblem}
  a(v, u_h) = L(v) \quad \forall v \in V_h^1,
\end{equation}
where $a$ and $L$ are bilinear and linear forms on $V_h^1 \times
V_h^2$ and $V_h^1$, respectively. Discretizing~\eqref{eq:varproblem},
one obtains a linear system $A U = b$ for the degrees of freedom~$U$
of the discrete solution~$u_h$.

In general, we shall be concerned with the discretization of a general
multilinear form of arity~$\rho$,
\begin{equation} \label{eq:a}
  a : V_h^1 \times V_h^2 \times \dots \times V_h^{\rho} \rightarrow \R.
\end{equation}
Typically, the arity is~$\rho = 1$ (linear forms) or $\rho = 2$
(bilinear forms), but forms of higher arity also appear
(see~\cite{KirLog2006}). For illustration purposes, we consider the
discretization of the mixed Poisson problem~\eqref{eq:mixedpoisson} in
the following example.
\begin{example}[discrete mixed Poisson]{\rm
  Let $\Sigma_h$ and $W_h$ be discrete spaces approximating $\HdivO$
  and $L^2(\Omega)$, respectively. We may then
  write~\eqref{eq:mixedpoisson} in the canonical
  form~\eqref{eq:varproblem} by defining
  \begin{subequations}
    \label{eq:discretemixedpoisson}
    \begin{align}
      a((\tau_h, v_h), (\flux_h, u_h)) &=
      \langle \tau_h, \flux_h \rangle - \langle \Div \tau_h, u_h \rangle
      + \langle v_h, \Div \flux_h \rangle, \\
      L((\tau_h, v_h)) &= \langle v_h, f \rangle
    \end{align}
  \end{subequations}
  for $ (\tau_h, v_h) \in V^1_h = \Sigma_h \times W_h$ and
  $ (\flux_h, u_h) \in V^2_h = V^1_h$.}
\end{example}

To discretize the multilinear form~\eqref{eq:a}, we let
$\{\phi^j_k\}_{k=1}^{N_j}$ denote a basis for $V^j_h$ for $j=1,
2, \ldots, \rho$ and we define the global tensor
\begin{equation}
  A_i = a(\phi_{i_1}^1, \phi_{i_2}^2, \dots,
  \phi_{i_{\rho}}^{\rho}),
\end{equation}
where $i = (i_1, i_2, \dots, i_{\rho})$ is a multi-index. Throughout,
$j, k$ denote simple indices. If the multilinear form is defined as an
integral over $\Omega = \cup_{K\in\mathcal{T}_h} K$, the tensor~$A$
may be computed by assembling the contributions from all elements,
\begin{align}
  A_{i} = a(\phi_{i_1}^1, \phi_{i_2}^2, \dots,
  \phi_{i_{\rho}}^{\rho})
  = \sum_{K \in \triang_h} a^K(\phi_{i_1}^1,
  \phi_{i_2}^2, \dots, \phi_{i_{\rho}}^{\rho}),
\end{align}
where $a^K$ denotes the contribution from element $K$.  We further let
$\{\phi_k^{K,j}\}_{k=1}^{n_j}$ denote the local finite element basis for
$V^j_h$ on $K$ and define the \emph{element tensor}~$A^K$ by
\begin{equation}
  A_i^K = a^K(\phi_{i_1}^{K, 1},
  \phi_{i_2}^{K, 2}, \dots, \phi_{i_{\rho}}^{K, \rho}).
\end{equation}
The assembly of the global tensor~$A$ thus reduces to the computation
of the element tensor $A^K$ on each element $K$ and the insertion of
the entries of $A^K$ into the global tensor $A$.

In~\cite{KirLog2007a}, it was shown that if the local basis on each
element~$K$ may be obtained as the image of a basis on a
\emph{reference element}~$K_0$ by the standard (affine) isomorphism
$\mathcal{F}_K: H^1(K_0) \rightarrow H^1(K)$, then the element tensor
$A^K$ may be represented as a tensor contraction of a \emph{reference
  tensor} $A^0$, depending only on the form $a$ and the reference
basis, and a \emph{geometry tensor} $G_K$, depending on the geometry
of the particular element~$K$,
\begin{equation} \label{eq:repr}
  A_i^{K} = A_{i\alpha}^0 G_K^{\alpha},
\end{equation}
with summation over the multi-index~$\alpha$. It was further
demonstrated in~\cite{KirLog2007a} that this representation may
significantly reduce the operation count for computing the element
tensor compared to standard evaluation schemes based on quadrature.

Below, we extend the representation~\eqref{eq:repr} to hold not only
for bases that may be affinely mapped from a reference element, but
also for finite element spaces that must be transformed by a Piola
mapping.

\subsection{A representation theorem}

We now state the general representation theorem for multilinear forms
on $H^1$, $\Hcurl$, $\Hdiv$ (and $L^2$). Instead of working out the
details of the proof here, we refer the reader to the proof presented
in~\cite{KirLog2007a} for $H^1$, and we illustrate the main points for
$\Hdiv$ and $\Hcurl$ by a series of examples.

\begin{thm}
  Let $K_0 \subset \R^n$ be a reference element and let $F_K: K_0
  \rightarrow K = F_K(K_0)$ be a nondegenerate, affine mapping with
  Jacobian $J_K$. For $j = 1, 2, \ldots, \rho$, let $\{\phi_{k}^{K, j}
  \}_{k}$ denote a basis on $K$ generated from a reference basis
  $\{\Phi_{k}^{j}\}_{k}$ on $K_0$, that is, $\phi_{k}^{K, j} =
  \mathcal{F}_{K}^j (\Phi_{k}^{j})$, where $\mathcal{F}_{K}^j$ is
  either of the mappings defined by~\eqref{eq:affine},
  \eqref{eq:contraPiola}, or~\eqref{eq:coPiola}.

  Then there exists a reference tensor $A_{i}^0$, independent of
  $K$, and a geometry tensor $G_K$ such that
  $  A^{K} = A^0 : G_K$,
  that is,
  \begin{equation}
    A^K_i = \sum_{\alpha\in\mathcal{A}} A^0_{i\alpha} G_K^{\alpha}
    \quad \forall i \in \mathcal{I},
  \end{equation}
  for a set of \emph{primary indices}~$\mathcal{I}$ and
  \emph{secondary indices}~$\mathcal{A}$.  In fact, the reference
  tensor $A^0$ takes the following canonical form:
  \begin{equation}
    A^0_{i\alpha}
    =
    \sum
    \int_{K_0}
    \prod_j
    D_X^{(\cdot)}
    \Phi^j_{(\cdot)}[(\cdot)]
    \dX;
  \end{equation}
  that is, it is the sum of integrals of products of basis function
  components and their derivatives on the reference element~$K_0$, and
  the geometry tensor $G_K$ is the outer product of the
  coefficients~$c_{(\cdot)}$ of any weight functions with a tensor
  that depends only on the Jacobian $J_K$,
  \begin{equation}
    G_K^{\alpha}
    =
    \prod
    c_{(\cdot)} \,
    \frac{|\det J_K|}{(\det J_K)^{\gamma}}
    \sum
    \prod
    \frac{\partial X_{(\cdot)}}{\partial x_{(\cdot)}}
    \prod
    \frac{\partial x_{(\cdot)}}{\partial X_{(\cdot)}} ,
  \end{equation}
  for some integer $\gamma$.
\end{thm}

\subsection{Examples}

To this end, we start by considering the vector-valued $L^2(\Omega)$
inner product, defining a bilinear form:
\begin{equation} \label{eq:mass}
  a(v, u) = \int_{\Omega} v \cdot u \dx.
\end{equation}
In the following, we let $x$ denote coordinates on~$K$ and we let $X$
denote coordinates on the reference element $K_0$. $F_K$ is an affine
mapping from $K_0$ to $K$, that is,\break $x = F_K(X) = J_K X + x_K$. We
further let $\phi^K$ denote a field on $K$ obtained as the image of a
field $\Phi$ on the reference element $K_0$,
$\phi^{K} = \mathcal{F}_K^{(\cdot)}(\Phi)$.
We aim to illustrate the differences and similarities of the
representations of the mass matrix for different choices of mappings
$\mathcal{F}_K$, in particular, affine, contravariant Piola, and
covariant Piola.
\begin{example}[the mass matrix with affinely mapped basis]{\rm
  Let $\mathcal{F}_K$ be the affine mapping,
  $  \mathcal{F}_K(\Phi) = \Phi \circ F_K^{-1}$.
  Then, the element matrix $A^K$ for~\eqref{eq:mass} is given by
  \begin{equation}
    A^K_i = \int_{K} \phi_{i_1}^{K, 1}(x) \cdot \phi_{i_2}^{K, 2}(x) \dx =
    |\det J_K|  \int_{K_0}
    \Phi^1_{i_1}[\beta](X) \, \Phi^2_{i_2}[\beta](X) \dX,
  \end{equation}
  where we use $\Phi[\beta]$ to denote component $\beta$ of the
  vector-valued function $\Phi$ and implicit summation over the
  index~$\beta$. We may thus represent the element matrix as the
  tensor contraction~\eqref{eq:repr} with reference and geometry
  tensors given by}
  \begin{align*}
    A^0_{i} &=
    \int_{K_0} \Phi^1_{i_1}[\beta](X) \, \Phi^2_{i_2}[\beta](X) \dX, \\
    \quad \quad
    G^K &= | \det J_K|.
  \end{align*}
\end{example}
\unskip%

We proceed to examine the representation of the mass matrix when the
basis functions are transformed with the contravariant and the
covariant Piola transforms.
\begin{example}[the mass matrix with contravariantly mapped basis]{\rm
  Let $\contraPiola_{K}$ be the contravariant Piola mapping,
  \begin{equation*}
    \contraPiola_{K}(\Phi) = \frac{1}{\det J_K} J_K \Phi \circ F_K^{-1}.
  \end{equation*}
  Then, the element matrix $A^K$ for~\eqref{eq:mass} is given by
  \begin{align*}
    A^K_i &= \int_{K} \phi_{i_1}^{K,1}(x) \cdot \phi_{i_2}^{K,2}(x) \dx \\
    &= \frac{|\det J_K|}{(\det J_K)^2} \frac{\px_{\beta}}{\pX_{\alpha_1}}
    \frac{\px_{\beta}}{\pX_{\alpha_2}}
    \int_{K_0} \Phi^1_{i_1}[\alpha_1](X) \, \Phi^2_{i_2}[\alpha_2](X) \dX.
  \end{align*}
  We may thus represent the element matrix as the
  tensor contraction~\eqref{eq:repr} with reference and geometry
  tensors given by}
  \begin{align*}
    A^0_{i\alpha} &= \int_{K_0} \Phi^1_{i_1}[\alpha_1](X) \,
    \Phi^2_{i_2}[\alpha_2](X) \dX, \\
    G_{\alpha}^K &= \frac{|\det J_K|}{(\det J_K)^2} \frac{\px_{\beta}}{\pX_{\alpha_1}}
    \frac{\px_{\beta}}{\pX_{\alpha_2}}.
  \end{align*}
\end{example}
\begin{example}[the mass matrix with covariantly mapped basis]{\rm
  Let $\coPiola_{K}$ be the covariant Piola mapping,
  \begin{equation*}
    \coPiola_{K}(\Phi) = J_K^{-T} \Phi \circ F_K^{-1}.
  \end{equation*}
  Then, the element tensor (matrix) $A^K$ for~\eqref{eq:mass} is given
  by
  \begin{align*}
    A^K_i &= \int_{K} \phi_{i_1}^{K,1}(x) \cdot \phi_{i_2}^{K,2}(x) \dx \\
    &=  |\det J_K| \, \frac{\pX_{\alpha_1}}{\px_{\beta}}
    \frac{\pX_{\alpha_2}}{\dx_{\beta}}
    \int_{K_0} \Phi^1_{i_1}[\alpha_1](X) \,
    \Phi^2_{i_2}[\alpha_2](X) \dX.
  \end{align*}
  We may thus represent the element matrix as the
  tensor contraction~\eqref{eq:repr} with reference and geometry
  tensors given by}
  \begin{align*}
    A^0_{i\alpha} &= \int_{K_0} \Phi^1_{i_1}[\alpha_1](X) \,
    \Phi^2_{i_2}[\alpha_2](X) \dX, \\
    G_{\alpha}^K &= |\det J_K| \, \frac{\pX_{\alpha_1}}{\px_{\beta}}
    \frac{\pX_{\alpha_2}}{\px_{\beta}}.
  \end{align*}
\end{example}
\unskip%

\looseness=1 We observe that the representation of the mass matrix differs for
affine, contravariant Piola, and covariant Piola. In particular, the
geometry tensor is different for each mapping, and the reference
tensor has rank two for the affine mapping, but rank four for the
Piola mappings. We also note that the reference tensor for the mass
matrix in the case of the covariant Piola mapping transforms in the
same way as the reference tensor for the stiffness matrix in the case
of an affine mapping (see~\cite{KirLog2007a}).

It is important to consider the storage requirements for this tensor
contraction approach and when other approaches might be appropriate.
For either the \( H(\mathrm{div}) \) or \( H(\mathrm{curl}) \) mass
matrix, for example, the reference tensor \( A^0 \) has rank four (two
indices for vector components and two for basis functions).  As such,
the storage requirements for \( A^0 \) are \( d^2 n^2 \), where \( d
=2,3 \) is the spatial dimension and \( n \) is the number of
reference element basis functions.  We also note that \( n =
\mathcal{O}(r^d) \), where \( r \) is the polynomial degree.  Storing
\( A^0 \) is thus comparable to storing \( d^2 \) element mass
matrices.  This is a modest, fixed amount of storage, independent of
the mesh. The tensor contraction may be computed in several different
ways.  The default option used by FFC is to generate straightline code
for performing the contraction of \( A^0 \) and \( G^K
\). Alternatively, one may also consider \( A^0 \) being stored in
memory as an array and applied via BLAS.  In the first case, the size
of generated code can become a problem for complex forms or high-order
methods, although this is not as large of a problem in the second
case.  The geometry tensor, \( G^K \), must be computed for each
element of the mesh.  For either the contravariant or covariant case,
\( G^K \) is a \( d \times d \) array and so is comparable to storing
the cell Jacobian for each cell of the mesh.  For more complicated
forms, storing \( G^K \) for each cell can become more expensive.
However, FFC currently stores only one such \( G^K \) at a time,
interleaving construction of \( G^K \) and its multiplication by \(
A^0 \).  For more complicated bilinear forms (such as ones involving
multiple material coefficients), the memory requirements of \( A^0 \)
and \( G^K \) both grow with the polynomial degree, which can lead to
inefficiency relative to a more traditional, quadrature-based
approach.  For a thorough study addressing some of these issues, we
refer the reader to~\cite{OelWel09}.

FFC is typically used to form
a global sparse matrix, but for high-degree elements, static
condensation or matrix-free approaches will be more appropriate.
This is a result of the large number of internal degrees of freedom
being stored in the sparse matrix and is an artifact of assembling a
global matrix rather than our tensor contraction formulation as such.

We conclude by demonstrating how the divergence term
from~\eqref{eq:discretemixedpoisson} is transformed with the
contravariant Piola (being the relevant mapping for $\Hdiv$).
\begin{example}[divergence term]
  \label{ex:div,simplification}{\rm
  Let $\mathcal{F}_K$ be the affine mapping, let $\contraPiola_K$ be
  the contravariant Piola mapping, and consider the bilinear form
  \begin{equation} \label{eq:divform}
    a(v, \sigma) = \int_K v \, \Div \sigma \dx
  \end{equation}
  for $(v, \sigma)  \in V^1 \times V^2$.  Then, if $\phi^{K, 1} =
  \mathcal{F}_K(\Phi^1)$ and $\phi^{K, 2} = \contraPiola_{K}(\Phi^2)$,
  the element matrix $A^K$ for~\eqref{eq:divform} is given by
  \begin{equation*}
    A^K_i = \int_{K} \phi^{K, 1}_{i_1} \,  \Div \phi^{K, 2}_{i_2} \dx
    = \frac{|\det J_K|}{\det J_K} \frac{\px_{\beta}}{\pX_{\alpha_1}}
    \frac{\pX_{\alpha_2}}{\px_{\beta}}
    \int_{K_0} \Phi^1_{i_1} \, \frac{\partial\Phi^2_{i_2}[\alpha_1]}{\partial X_{\alpha_2}} \dX.
  \end{equation*}
  Noting that $\frac{\px_{\beta}}{\pX_{\alpha_1}}
  \frac{\pX_{\alpha_2}}{\px_{\beta}} = \delta_{\alpha_1\alpha_2}$,
  we may simplify to obtain
  \begin{equation*}
    A^K_i = \frac{|\det J_K|}{\det J_K}
    \int_{K_0} \Phi^1_{i_1} \, \frac{\partial\Phi^2_{i_2}[\alpha_1]}{\partial
    X_{\alpha_1}} \dX
    = \pm
    \int_{K_0} \Phi^1_{i_1} \, \Div \Phi^2_{i_2} \dX.
  \end{equation*}
  We may thus represent the element matrix as the
  tensor contraction~\eqref{eq:repr} with reference and geometry
  tensors given by}
  \begin{align*}
    A^0_i &= \int_{K_0} \Phi^1_{i_1} \, \Div \Phi^2_{i_2} \dX, \\
    G_{\alpha}^K &= \pm 1.
  \end{align*}
\end{example}
\unskip%

The simplification in the final example is a result of the
isomorphism, induced by the contravariant Piola transform, between
$H(\Div, K_0)$ and $H(\Div, K)$. FFC takes special
care of such and similar simplifications.

\section{Assembling \boldmath$\Hdiv$ and $\Hcurl$ elements}
\label{sec:assembling}

To guarantee global continuity with Piola-mapped elements, special
care has to be taken with regard to the numbering and orientation of
geometric entities, in particular the interplay between local and
global orientation. This is well known, but is rarely discussed in the
standard references, though some details may be found in~\cite{Mon03,
Sch03}. We discuss here some of these issues and give a strategy for
dealing with directions of normals and tangents that simplifies
assembly over $\Hdiv$ and $\Hcurl$. In fact, we demonstrate that one
may completely remove the need for contending with directions by using
an appropriate numbering scheme for the simplicial mesh.

\subsection{Numbering scheme}
\label{subsec:numbering}
The numbering and orientation of geometric entities in FFC follows the
UFC specification~\cite{AlnLog08}. In short, the numbering scheme
works as follows. A global index is assigned to each vertex of the
tessellation $\mathcal{T}_h$ (consisting of triangles or
tetrahedra). If an edge adjoins two vertices $v_i$ and $v_j$, we
define the direction of the edge as going from vertex $v_i$ to vertex
$v_j$ if $i < j$.  This gives a unique orientation of each edge.  The
same convention is used locally to define the directions of the local
edges on each element. Thus, if an edge adjoins the first and second
vertices of a tetrahedron, then the direction is from the first to the
second vertex. A similar numbering strategy is employed for faces. The
key now is to require that the vertices of each element are always
ordered based on the their global indices.

\begin{figure}
  \begin{center}
{\hspace*{-3.5pt}{    \psfrag{0}{$1$}
    \psfrag{1}{$2$}
    \psfrag{2}{$3$}
    \psfrag{K0}{$K_0$}
    \psfrag{K1}{$K$}
    \psfrag{K2}{$K'$}
    \psfrag{F1}{$F_K$}
    \psfrag{F2}{$F_{K'}$}
    \psfrag{e}{$e$}
    \psfrag{E1}{$E$}
    \psfrag{E2}{$E'$}
    \psfrag{v10}{$v_{10}$}
    \psfrag{v20}{$v_{20}$}
    \psfrag{v50}{$v_{50}$}
    \psfrag{v75}{$v_{75}$}
    \includegraphics[width=\textwidth]{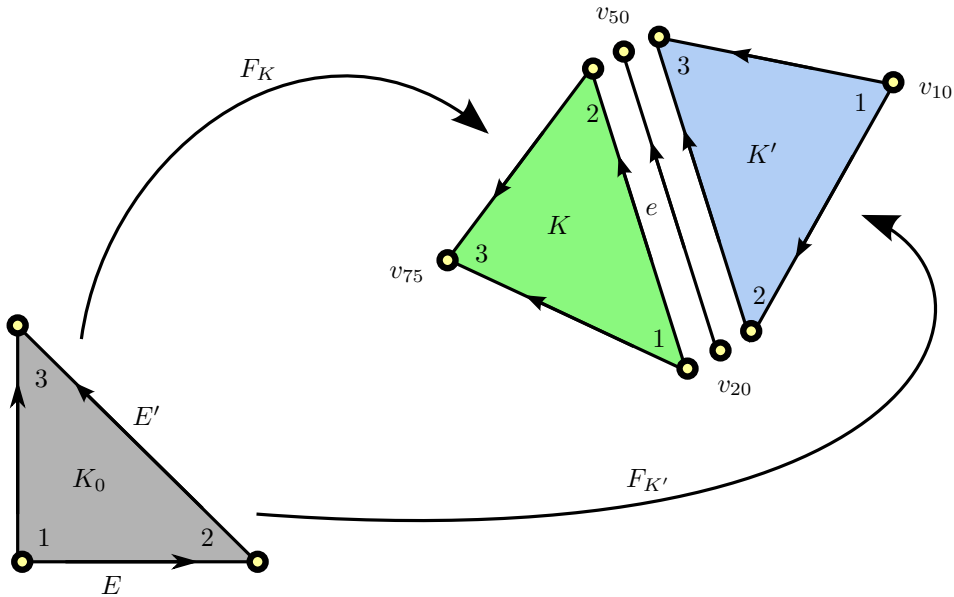}}}
    \caption{Two adjacent triangles will always agree on the direction
      of a common edge tangent or normal. The two triangles in the
      figure share a common edge between the global vertices $v_{20}$
      and $v_{50}$. These two vertices have different local indices
      for $K$ ($1, 2$) and $K'$ ($2, 3$), but the ordering convention,
      local numbering according to ascending global indices, ensures
      that both triangles agree on the direction of the common edge
      $e$.}
    \label{fig:patchofelements}
  \end{center}
\end{figure}

For illustration, consider first the two-dimensional case. Let $K_0$
be the UFC reference triangle, that is, the triangle defined by the
vertices $\{(0,0), (1,0), (0, 1)\}$. Assume that $K= F_K(K_0)$ and $K'
= F_{K'}(K_0)$ are two physical triangles sharing an edge $e$ with
normal $n$. If $e$ adjoins vertices $v_i$ and $v_j$ and is directed
from $v_i$ to $v_j$, it follows from the numbering scheme that $i <
j$. Since the vertices of both $K$ and $K'$ are ordered based on their
global indices, and the local direction (as seen from $K$ or $K'$) of
an edge is based on the local indices of the vertices adjoining that
edge, this means that the local direction of the edge $e$ will agree
with the global direction, both for $K$ and $K'$.  Furthermore, if we
define edge normals as clockwise rotated tangents, $K$ and $K'$ will
agree on the direction of the normal of the common edge. The reader is
encouraged to consult Figure~\ref{fig:patchofelements} for an
illustration.

\begin{figure}
\vspace*{-8pt}
\begin{center}
{\hspace*{0.75pc}{%
    \psfrag{1}{$1$}
    \psfrag{2}{$2$}
    \psfrag{3}{$3$}
    \psfrag{4}{$4$}
    \psfrag{v10}{$v_{10}$}
    \psfrag{v20}{$v_{20}$}
    \psfrag{v50}{$v_{50}$}
    \psfrag{v75}{$v_{75}$}
    \psfrag{v92}{$v_{92}$}
    \psfrag{e1}{$e_1$}
    \psfrag{e2}{$e_2$}
    \psfrag{ee1}{$e_1'$}
    \psfrag{ee2}{$e_2'$}
    \psfrag{K1}{$K$}
    \psfrag{K2}{$K'$}
    \includegraphics[width=\textwidth]{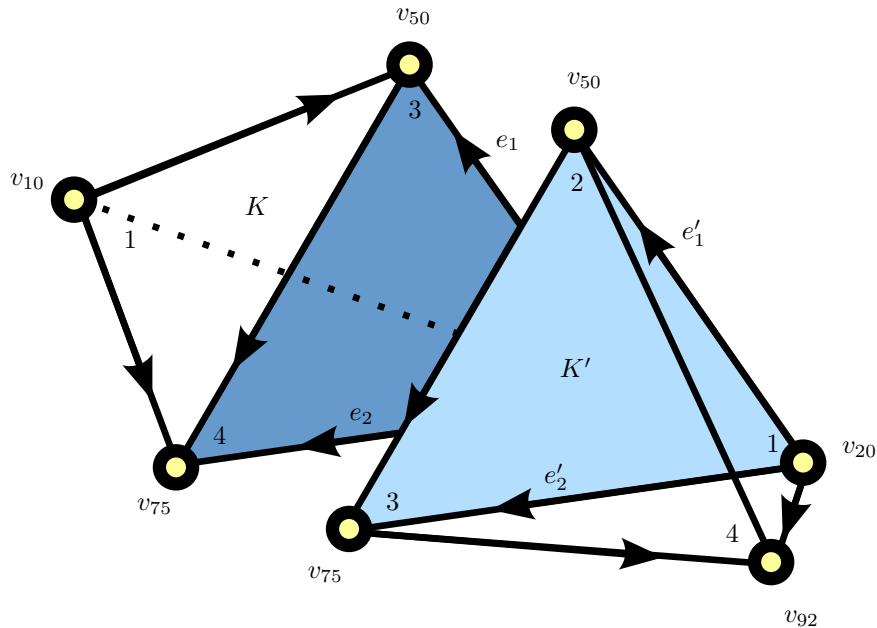}}}
    \caption{Two adjacent tetrahedra will always agree on the
      direction of a common edge tangent or face normal.  The two
      tetrahedra in the figure share a common face defined by the
      global vertices $v_{20}$, $v_{50}$, and $v_{75}$. These three
      vertices have different local indices for $K$ ($2, 3, 4$) and $K'$
      ($1, 2, 3$), but the ordering convention, local numbering according
      to ascending global indices, ensures that both triangles agree
      on the direction of the common edges. In particular, the two
      tetrahedra agree on the directions of the first two edges of the
      common face and the direction of the face normal $n \propto e_1
      \times e_2 = e_1' \times e_2'$.}
    \label{fig:patchofelements:3D}
  \end{center}
\vspace*{-6pt}
\end{figure}

The same argument holds for the direction of edges and face normals in
three dimensions. In particular, if face normals are consistently
defined in terms of edges, it is straightforward to ensure a common
direction. Consider two tetrahedra $K$ and $K'$ sharing a face $f$,
defined by three vertices $v_{i_1}, v_{i_2}, v_{i_3}$ such that $i_1 <
i_2 < i_3$.  Clearly, $v_{i_1}$ will be the vertex with the lowest
index of the face $f$ for both $K$ and $K'$. Furthermore, each of the
two edges that adjoin $v_{i_1}$, that is, the edge from $v_{i_1}$ to
$v_{i_2}$ and the edge from $v_{i_1}$ to $v_{i_3}$, has a unique
direction by the previous arguments. These two edges can therefore
define consistent tangential directions of the face. Taking the
normalized cross-product of these edges gives a consistent face
normal. This is the approach used by FIAT/FFC. As a consequence, two
adjacent tetrahedra sharing a common face will always agree on the
direction of the tangential and normal directions of that face. This
is illustrated in Figure~\ref{fig:patchofelements:3D}.

We emphasize that the numbering scheme above does not result in a
consistent orientation of the boundary of each element. It does,
however, ensure that two adjacent elements sharing a common edge or
face will always agree on the orientation of that edge or face. In
addition to facilitating the treatment of tangential and normal
traces, a unique orientation of edges and faces simplifies assembly of
higher order Lagrange elements. A similar numbering scheme is proposed
in the monograph~\cite{Mon03} for tetrahedra in connection with
$\Hcurl$ finite elements. Also, we note that the numbering scheme and
the consistent facet orientation that follows render only one
reference element necessary, in contrast to the approach
of~\cite{art:AinCoy2003}.

\subsection{Mapping nodal basis functions}
Next, we show how this numbering scheme and the FIAT choice of degrees
of freedom give the necessary $\Hdiv$ or $\Hcurl$ continuity. Assume
that we have defined a set of nodal basis functions on $K_0$, that is,
$\{\Phi_i\}_{i=1}^n$ such that
\begin{equation*}
  \ell_i(\Phi_j) = \delta_{ij}, \quad i,j=1,2,\ldots,n,
\end{equation*}
for a set of degrees of freedom $\{\ell_i\}_{i=1}^n$. These basis
functions are mapped to two physical elements $K$ and $K'$ by an
appropriate transformation $\mathcal{F}$ (contravariant or covariant
Piola), giving a set of functions on $K$ and $K'$, respectively. We
demonstrate below that as a consequence of the above numbering scheme,
these functions will indeed be the restrictions to $K$ and $K'$ of an
appropriate global nodal basis.

Consider $\Hcurl$ and a global degree of freedom~$\ell$ defined
as the tangential component at a point~$x$ on a global edge~$e$ with
tangent~$t$, weighted by the length of the edge~$e$,
\begin{equation*}
  \ell(v) = \|e\| \, v(x) \cdot t = v(x) \cdot e.
\end{equation*}
Let $\coPiola$ be the covariant Piola mapping as before and let
$\phi^{K}$ and $\phi^{K'}$ be two basis functions on $K$ and $K'$
obtained as the mappings of two nodal basis functions, say, $\Phi$ and
$\Phi'$, on $K_0$,
\begin{equation*}
  \phi^K = \coPiola_K(\Phi) \quad \text{ and } \quad \phi^{K'} = \coPiola_{K'}(\Phi').
\end{equation*}
Assume further that $\Phi$ is the nodal basis function corresponding
to evaluation of the tangential component at the point $X\in K_0$
along the edge $E$, and that $\Phi'$ is the nodal basis function
corresponding to evaluation of the tangential component at the point
$X'\in K_0$ along the edge $E'$. Then, if $x = F_K(X) = F_{K'}(X')$,
the covariant Piola mapping ensures that
\begin{equation*}
  \phi^K(x) \cdot e = \Phi(X) \cdot E = 1
  \quad \text{ and } \quad
  \phi^{K'}(x') \cdot e' = \Phi'(X') \cdot E' = 1.
\end{equation*}
Thus, since $e = e'$, it follows that
\begin{equation*}
  \ell(\phi^K) = \ell(\phi^{K'}).
\end{equation*}
Continuity for $\Hdiv$ may be demonstrated similarly.

In general, FFC allows elements for which the nodal basis on the
reference element $K_0$ is mapped exactly to the nodal basis for each
element $K$ under some mapping $\mathcal{F}$, whether this be affine
change of coordinates or one of the Piola transformations.  While this
enables a considerable range of elements, as considered in this paper,
it leaves out many other elements of interest.  As an example, the
Hermite triangle or tetrahedron~\cite{Cia76} does not transform
equivalently.  The Hermite triangle has degrees of freedom which are
point values at the vertices and the barycenter, and the partial
derivatives at each vertex.  Mapping the basis function associated
with a vertex point value affinely yields the correct basis function
for \( K \), but not for the derivative basis functions.  A simple
calculation shows that a function with unit \( x \)-derivative and
vanishing \( y \)-derivative at a point generally maps to a function
for which this is not the case.  In fact,  the
function value basis functions transform affinely, but the pairs of
derivative basis functions at each vertex must be transformed
together; that is, a linear combination of their image yields the correct basis
functions.

Examples of other elements requiring more general types of mappings include
the scalar-valued
Argyris and Morley elements as well as the Arnold--Winther symmetric
elasticity element~\cite{ArnWin02} and the Mardal--Tai--Winther element
for Darcy--Stokes flow~\cite{MarTai02}.  Recently, a special-purpose
mapping for the Argyris element has been developed by Dom\'{\i}nguez and
Sayas~\cite{DomSay08}, and we are generalizing this work as an
extension of the FIAT project as outlined below.

If $\{ \Phi_i \}$ is the reference finite element basis and $\{
\phi^K_i \}_i$ is the physical finite element basis, then equivalent
elements satisfy $\phi^K_i = \mathcal{F}_K(\Phi_i)$ for each~$i$. If
the elements are not equivalent under \( \mathcal{F}_K \), then \( \{
\phi^K_i \}_i \) and \( \{ \mathcal{F}_K(\Phi_i) \}_{i} \) form two
different bases for the polynomial space. Consequently, there exists a
matrix $M^K$ such that $\phi^K_i = \sum_j M^K_{ij}
\mathcal{F}_K(\Phi_j)$.  In the future, we hope to extend FIAT to
construct this matrix $M$ and FFC to make use of it in constructing
variational forms, further extending the range of elements available
to users.

\subsection{A note about directions}
An alternative orientation of shared facets gives rise to a special
case of such transformations. It is customary to direct edges in a
fashion that gives a consistent orientation of the boundary of each
triangle. However, this would mean that two adjacent triangles may
disagree on the direction of their common edge. In this setting,
normals would naturally be directed outward from each triangle, which
again would imply that two adjacent triangles disagree on the
direction of the normal on a common edge. It can be demonstrated that
it is then more appropriate to define the contravariant Piola mapping
in the following slightly modified form:
\begin{equation*}
\contraPiola(\Phi) = \frac{1}{|\det J_K|} J_K \Phi \circ F_K^{-1};
\end{equation*}
that is, the determinant of the Jacobian appears without a sign.

To ensure global continuity, one would then need to introduce
appropriate sign changes for the mapped basis functions.  For two
corresponding basis functions $\phi^K$ and $\phi^{K'}$ as above, one
would change the sign of $\phi^{K'}$ or $\phi^K$ such that both basis
functions correspond to the same global degree of freedom. Thus, one
may consider obtaining the basis functions on the physical element by
first mapping the nodal basis functions from the reference element
and then correcting those basis functions with a change of sign:
\begin{align*}
  \tilde{\phi}^K &= \mathcal{F}(\Phi), \\
  \phi^K         &= \pm \, \tilde{\phi}^K.
\end{align*}
This would correspond to a diagonal \( M^K \) transformation where the
entries are all \( \pm 1 \).

Since a multilinear form is linear in each of its arguments, this
approach corresponds to first computing a tentative element
tensor~$\tilde{A}^K$ and then obtaining $A^K$ from $\tilde{A}^K$ by a
series of rank one transforms. However, this procedure is unnecessary
if the contravariant Piola mapping is defined as
in~\eqref{eq:contraPiola} and the numbering scheme described in
section \ref{subsec:numbering} is employed.

For nonsimplicial meshes, such as meshes consisting of quadrilaterals
or hexahedra, the situation is somewhat more complicated. It is not
clear how to ensure a consistent, common local and global direction
for the edges. Therefore, the UFC specification instead requires a
consistent orientation of the boundary of each cell. In this
situation, the alternative approach, relying on the introduction of
sign changes, is more appropriate.

\section{Examples}
\label{sec:examples}

In order to demonstrate the veracity of the implementation and the
ease with which the $\Hdiv$ and $\Hcurl$ conforming elements can be
employed, we now present a set of numerical examples and include the
FFC code used to define the variational forms. In particular, we
return to the examples introduced in section~\ref{sec:hdivcurl}, which
include formulations of the Hodge Laplace equations, the cavity
resonator eigenvalue problem, and the weak symmetry formulation for
linear elasticity.

\subsection{The Hodge Laplacian}

Consider the weak formulations of the Hodge Laplace equation
introduced in Examples \ref{ex:mixedpoisson} and
\ref{ex:hodgelaplace}. For $\Omega \subset \R^2$ and differential $1$-
and $2$-forms, we have the mixed Poisson equation
\eqref{eq:mixedpoisson}. Stable choices of conforming finite element
spaces $\Sigma_h \times V_h \subset \Hdiv \times L^2$ include $V_h =
\DG_{r-1}$ in combination with $\Sigma_h \in \{ \RT_{r-1}, \BDFM_{r},
\BDM_{r} \}$ for $r = 1, 2, \dots$. The FFC code corresponding to the
latter choice of elements is given in Table
\ref{tab:codemixedpoisson}. Further, for $\Omega \subset \R^3$, we
give the FFC code for the formulation of \Eq \eqref{eq:curldiv} with
the element spaces $\NED_{r-1}^1 \times \RT_{r-1} \subset \Hcurl
\times \Hdiv$ in Table \ref{tab:codecurldiv}.

For testing purposes, we consider a regular tessellation of the unit
square/cube, $\Omega = [0, 1]^n$, $n = 2, 3$, and a given smooth
source for the two formulations. In particular, for
\eqref{eq:mixedpoisson}, we solve for
\begin{equation}
  u(x_1, x_2) = C \sin(\pi x_1) \sin(\pi x_2),
\end{equation}
with $C$ a suitable scaling factor, and for \eqref{eq:curldiv}, we let
\begin{equation}
  \label{eq:exactu:curldiv}
  u(x_1, x_2, x_3) =
  \begin{pmatrix}
    x_1^2 (x_1-1)^2 \sin(\pi x_2) \sin(\pi x_3) \\
    x_2^2 (x_2-1)^2 \sin(\pi x_1) \sin(\pi x_3) \\
    x_3^2 (x_3-1)^2 \sin(\pi x_1) \sin(\pi x_2)
  \end{pmatrix}.
\end{equation}
Note that $u$ given by \Eq \eqref{eq:exactu:curldiv} is
divergence-free and such that $u \times n = 0$ on the exterior
boundary, and thus satisfies the implicit natural boundary conditions
of \eqref{eq:curldiv}.

\begin{table}
\footnotesize
\caption{FFC code for the mixed Poisson equation.}  \vskip-5pt
  \label{tab:codemixedpoisson}
  \begin{code}{0.8}
    r = 3
    S = FiniteElement("BDM", "triangle", r)
    V = FiniteElement("DG", "triangle", r - 1)
    element = S + V

    (tau, v) = TestFunctions(element)
    (sigma, u) = TrialFunctions(element)

    a = (dot(tau, sigma) - dot(div(tau), u) + dot(v, div(sigma))*dx
    L = dot(v, f)*dx
  \end{code}
\end{table}
\begin{table}
\footnotesize
  \caption{FFC code for the {\rm{curl-div}} formulation of the Hodge Laplace equation.} \vskip-5pt
  \label{tab:codecurldiv}
  \begin{code}{0.8}
    r = 2
    CURL = FiniteElement("Nedelec", "tetrahedron", r - 1)
    DIV = FiniteElement("RT", "tetrahedron", r - 1)
    element = CURL + DIV

    (tau, v) = TestFunctions(element)
    (sigma, u) = TrialFunctions(element)

    a = (dot(tau, sigma) - dot(curl(tau), u) + dot(v, curl(sigma)) \
        + dot(div(v), div(u)))*dx
    L = dot(v, f)*dx
  \end{code}
\end{table}

\begin{figure}[!t]
  \begin{center}
    \includegraphics[width=0.45\textwidth]{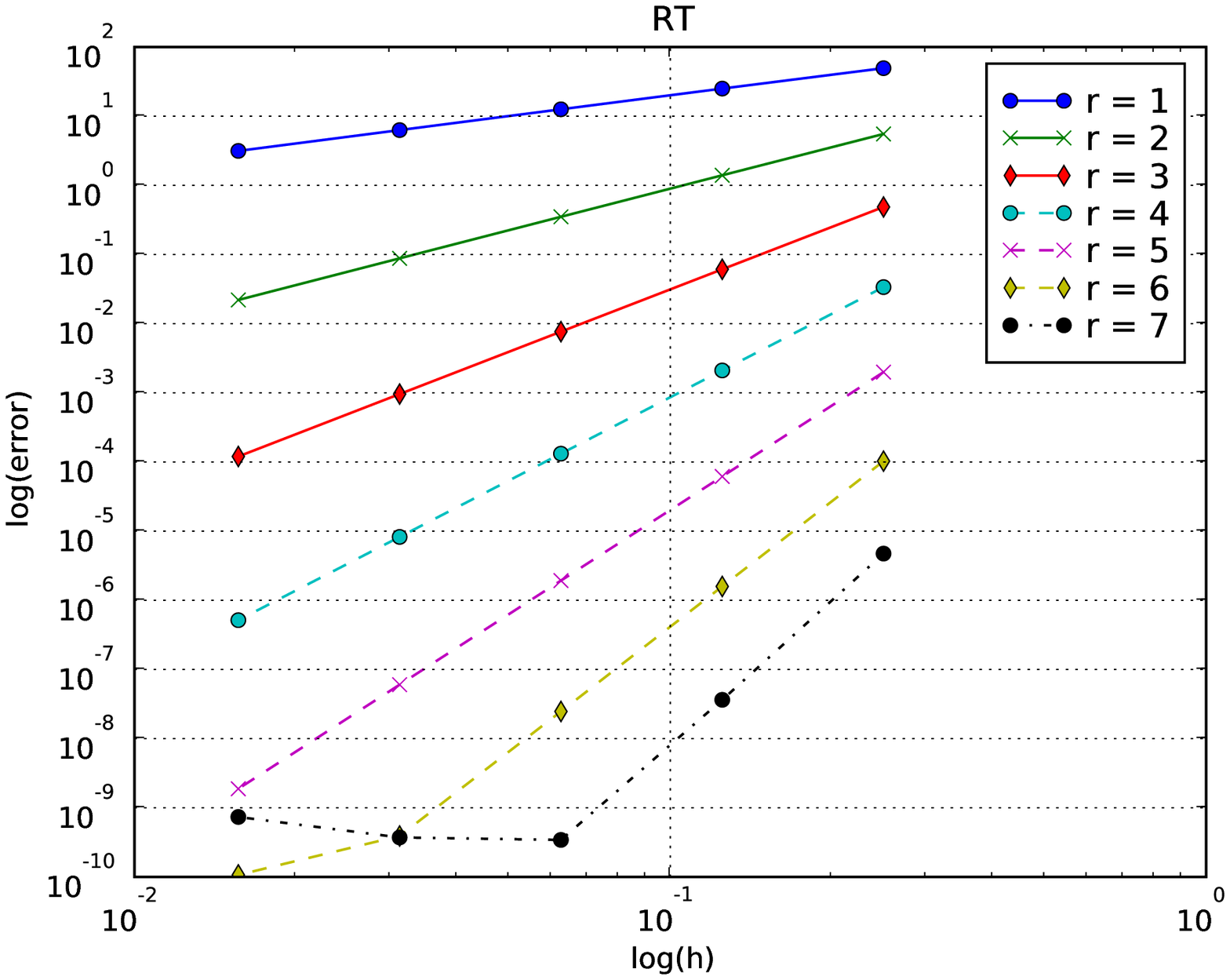}
    \includegraphics[width=0.45\textwidth]{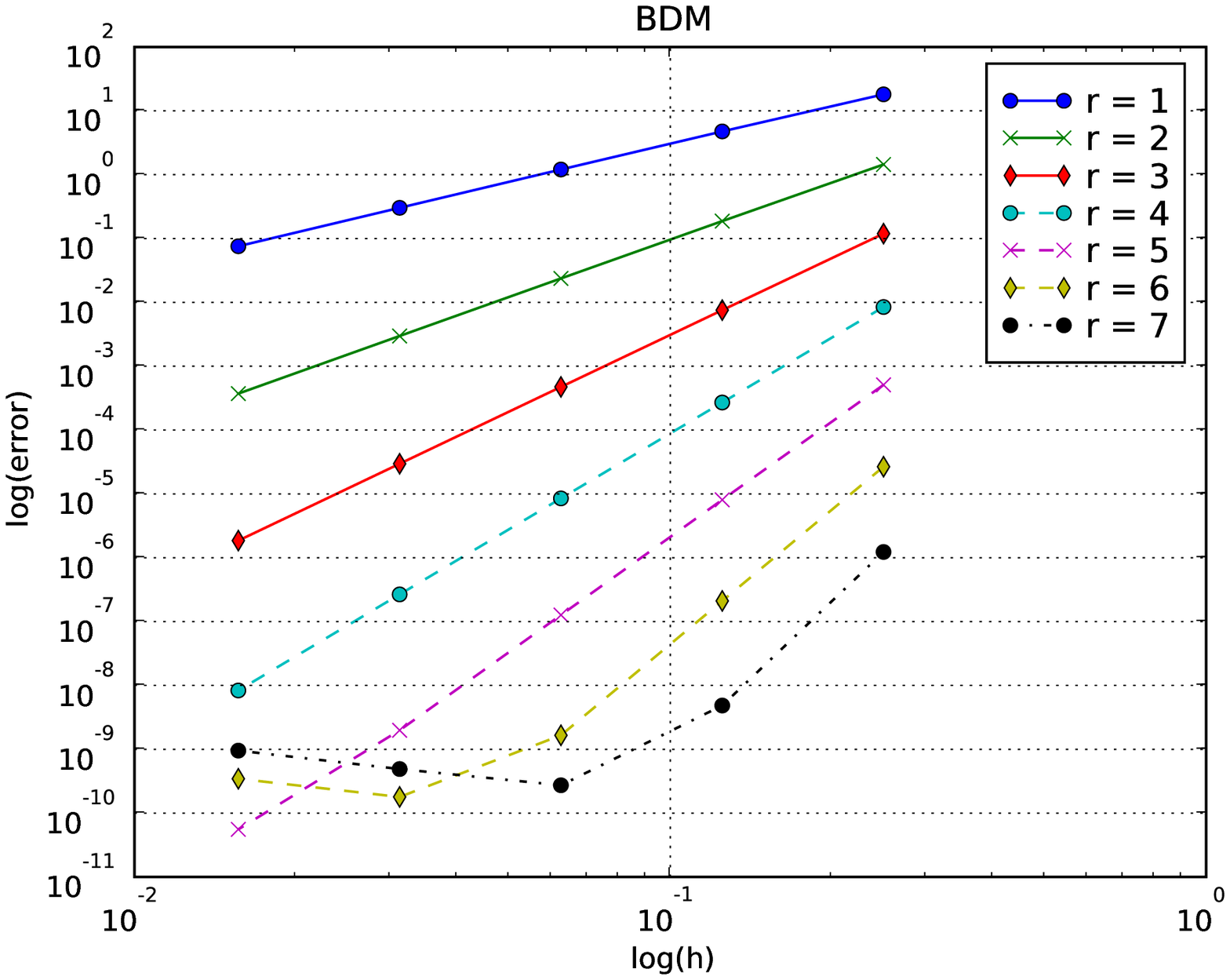}
    \caption{Convergence rates for the discretized mixed Poisson
    equation \Eq {\rm \eqref{eq:mixedpoisson}} using $\RT_{r-1} \times
    \DG_{r-1}$ (left) and $\BDM_{r} \times DG_{r-1}$ (right), $r = 1,2,\ldots,7$. Logarithmic plots of the $L^2$ error of the flux
    approximation: $||\sigma - \sigma_h||_0$ versus mesh size. The
    convergence rates in the left plot are $\mathcal{O}(h^r)$ and the
    convergence rates in the right plot are $\mathcal{O}(h^{r+1})$;
    cf.~Table {\rm \ref{tab:approximability}}. The error does not converge
    below $\sim10^{-10}$ in our experiments as a result of limited
    precision in the evaluation of integrals and/or linear
    solvers. The exact source of the limited precision has not been
    investigated in detail.}
    \label{fig:convergence_mxpoisson}
  \end{center}
\end{figure}

A comparison of the exact and the approximate solutions for a set of
uniformly refined meshes gives convergence rates in perfect agreement
with the theoretical values indicated by Table
\ref{tab:approximability}, up to a precision limit. Logarithmic plots
of the $L^2$ error of the flux using $\Sigma_h \in \{\RT_{r-1}, \BDM_r
\}$ versus the mesh size for $r = 1, 2, \dots, 7$ can be inspected in
Figure \ref{fig:convergence_mxpoisson} for the mixed Poisson problem
(with $C = 100$).

For the curl-div formulation of the Hodge Laplace equation \Eq
\eqref{eq:curldiv}, we have included convergence rates for $u$ and
$\sigma$ in Table \ref{table:convergence_curldiv}. Note that the
convergence rates for the combinations $\NED_{r} \times \RT_{r}$ and
$\NED_{r} \times \BDM_{r}$, $r = 1, 2$, are of the same order, except
for the $||\cdot||_{\Div}$ error of $u$, though the former combination
is computationally more expensive.

\subsection{The cavity resonator}

\begin{figure}[!t]
  \begin{center}
    \includegraphics[width=0.78\textwidth]{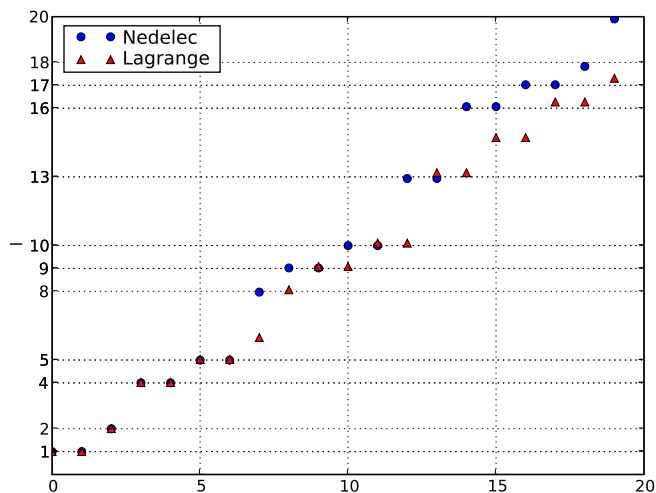}
    \caption{The first $20$ eigenvalues of the cavity resonator problem
      computed using first order \Nedelec{} elements ($\mathrm{NED}_0$)
      and Lagrange elements ($P_1$) on a coarse ($16\times 16$)
      criss-cross mesh. The exact analytical values are indicated by
      the horizontal grid lines.}
    \label{fig:eigenvalues}
  \end{center}
\end{figure}

{The analytical nonzero eigenvalues of the Maxwell eigenvalue problem
\eqref{eq:maxwell} with $\Omega = [0, \pi]^n$, $n = 2, 3$, are given by
\begin{equation}
\omega^2 = m_1^2 + m_2^2 + \cdots + m_n^2,\quad m_i \in \{ 0 \} \cup \mathbb{N},
\end{equation}
where at least $n - 1$ of the terms $m_i$ must be nonzero.  It is
well known~\cite{Mon03} that discretizations of this eigenvalue
problem using $H^1$ conforming finite elements produce spurious and
highly mesh-dependent eigenvalues $\omega^2$. The edge elements of the
\Nedelec{} type, however, give convergent approximations of the
eigenvalues. This phenomenon is illustrated in Figure
\ref{fig:eigenvalues}. There, the first $20$ nonzero eigenvalues
$\omega^2_{h,N}$ produced {\hfilneg}}

\clearpage

\begin{table}[!t]
\footnotesize
\caption{Averaged convergence rates for the discretized {\rm{curl-div}}
    formulation of the Hodge Laplace equation \Eq {\rm \eqref{eq:curldiv}}
    using $\NED_{r-1} \times \RT_{r-1}$, $r = 1, 2, 3$, and $\NED_{r}
    \times \BDM_{r}$, $r = 1, 2$. Number of degrees of freedom in the
    range $80{,}000$--$300{,}000$.}
    \label{table:convergence_curldiv}
  \centering
    \begin{tabular}{|l|c|c|c|c|}
      \hline
      \null\quad\, Element & $||\sigma-\sigma_h||_0$
      & $||\sigma-\sigma_h||_{\Curl}$ & $||u-u_h||_0$ & $||u-u_h||_{\Div}$ \\
      \hline
      $\NED_0 \times \RT_0$ & 0.99 & 0.98 & 0.99 & 0.98 \\
      $\NED_1 \times \BDM_1$ & 1.96 & 2.00 & 1.95 & 0.96 \\
      $\NED_1 \times \RT_1$ & 1.97 & 1.97 & 1.98 & 1.98 \\
      $\NED_2 \times \BDM_2$ & 3.00 & 2.99 & 2.97 & 1.97 \\
      $\NED_2 \times \RT_2$ & 2.98 & 2.96 & 2.97 & 2.97 \\
      \hline
    \end{tabular}
\end{table}

\clearpage

\noindent by the \Nedelec{} edge elements on a regular
criss-cross triangulation are given in comparison with the
corresponding Lagrange eigenvalue approximations
$\omega^2_{h,L}$. Note the treacherous spurious Lagrange
approximations such as $\omega^2_{h,L} \approx 6, 15$.

\subsection{Elasticity with weakly imposed symmetry}

As a final example, we consider a mixed finite element formulation of
the equations of linear elasticity with the symmetry of the stress
tensor imposed weakly as given in Example \ref{ex:elasticity}. In the
homogeneous, isotropic case, the inner product induced by the
compliance tensor $A$ reduces to
\begin{equation*}
  \langle \tau, A \sigma \rangle =
  \nu \langle \tau, \sigma \rangle - \zeta \langle \tr \tau, \tr \sigma \rangle
\end{equation*}
for $\nu, \zeta$ material parameters. A stable family of finite
element spaces for the discretization of \eqref{eq:elasticity} is
given by~\cite{ArnFalWin07}: $\BDM^2_r \times \DG^2_{r-1} \times
\DG_{r-1} \subset H(\Div; \Omega, \M) \times L^2(\Omega, \R^n)
\times L^2(\Omega)$, $r = 1, 2, \dots$. The $17$ lines of FFC code
sufficient to define this discretization are included in
Table~\ref{tab:codeelasticity}.

\begin{figure}[!t]
  \begin{center}
    \includegraphics[width=0.42\textwidth]{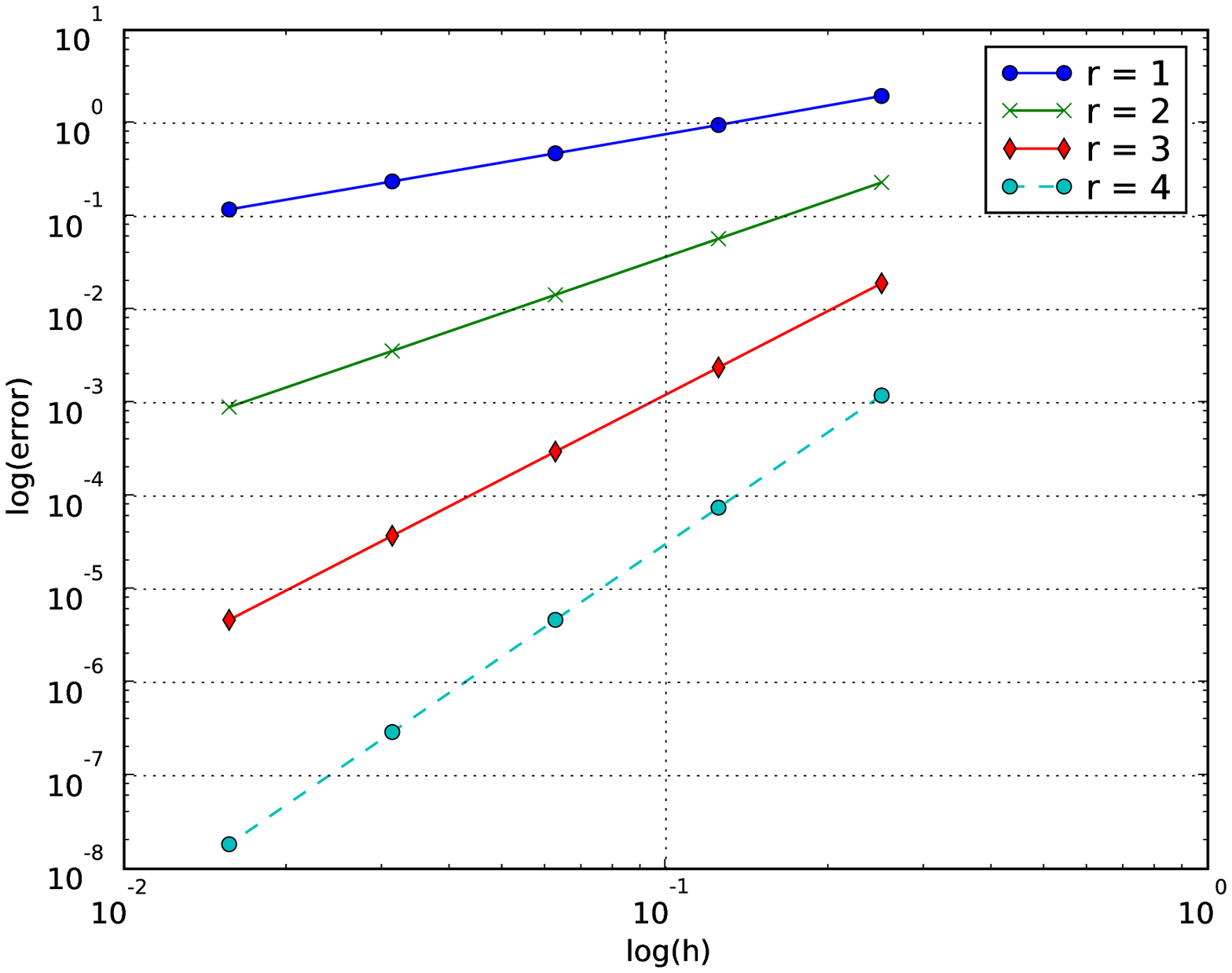}
    \includegraphics[width=0.43\textwidth]{eps/twin_dolphins.epsi}
    \caption{Left: Convergence rates for elastic stress approximations
      of {\rm \eqref{eq:elasticity}}. Logarithmic plot of $\Hdiv$ error of
      the approximated stress $\sigma$ versus mesh size. The
      convergence rates are $\mathcal{O}(h^r)$, $r = 1, 2, 3,
      4$. Right: Elastic dolphin hanging by the tail under a
      gravitational force.}
    \label{fig:elasticity}
  \end{center}
\vspace*{-8pt}
\end{figure}

Again to demonstrate convergence, we consider a regular triangulation
of the unit square and solve for the smooth solution
\begin{equation}
  u(x_0, x_1) =
  \begin{pmatrix}
    - x_1 \sin(\pi x_0) \\
    0.5 \pi x_1^2 \cos(\pi x_0)
  \end{pmatrix}.
\end{equation}
The theoretically predicted convergence rate of the discretization
introduced above is of the order $\mathcal{O}(h^r)$ for all computed
quantities. The numerical experiments corroborate this prediction. In
particular, the convergence of the stress approximation in the $\Hdiv$
norm can be examined in Figure~\ref{fig:elasticity}.
\begin{table}[!t]
\footnotesize \vskip5pt
\caption{FFC code for linear elasticity with weak symmetry.}
  \label{tab:codeelasticity}
  \begin{code}{0.8}
    def A(sigma, tau, nu, zeta):
        return (nu*dot(sigma, tau) - zeta*trace(sigma)*trace(tau))*dx

    def b(tau, w, eta):
        return (div(tau[0])*w[0] + div(tau[1])*w[1] + skew(tau)*eta)*dx

    nu =  0.5
    zeta =  0.2475
    r = 2

    S = FiniteElement("BDM", "triangle", r)
    V = VectorElement("Discontinuous Lagrange", "triangle", r-1)
    Q = FiniteElement("Discontinuous Lagrange", "triangle", r-1)
    MX = MixedElement([S, S, V, Q])

    (tau0, tau1, v, eta) = TestFunctions(MX)
    (sigma0, sigma1, u, gamma) = TrialFunctions(MX)
    sigma = [sigma0, sigma1]
    tau = [tau0, tau1]

    a = A(sigma, tau, nu, zeta) + b(tau, u, gamma) + b(sigma, v, eta)
    L = dot(v, f)*dx
  \end{code}
\end{table}

\section{Conclusions}
\label{sec:conclusion}

The relative scarcity of $\Hdiv$ and $\Hcurl$ mixed finite element
formulations in practical use may be attributed to their higher
theoretical and implementational threshold. Indeed, more care is
required to implement their finite element basis functions than the
standard Lagrange bases, and assembly poses additional difficulties.
However, as demonstrated in this work, the implementation of mixed
finite element formulations over $\Hdiv$ and $\Hcurl$ may be automated
and thus be used with the same ease as standard formulations over
$H^1$. In particular, the additional challenges in the assembly can be
viewed as not essentially different from those encountered when
assembling higher order Lagrange elements.

The efficiency of the approach has been further investigated by \O
lgaard and Wells~\cite{OelWel09}, with particular emphasis on the
performance when applied to more complicated PDEs. They conclude that
the tensor representation significantly improves performance for forms
below a certain complexity level, corroborating the previous results
of~\cite{KirLog2007a}. However, an automated, optimized quadrature
approach, also supported by FFC, may prove significantly better for
more complex forms. These findings indicate that a system for
automatically detecting the better approach may be valuable.

The tools (FFC, FIAT, DOLFIN) used to compute the results presented
here are freely available as part of the FEniCS
project~\cite{logg:www:03} and it is our hope that this may contribute
to further the use of mixed formulations in applications.


\vspace*{-4pt}
\begin{thebibliography}{10}

\bibitem{www:FEAP}
{\em {FEAP:} {A} {F}inite {E}lement {A}nalysis {P}rogram}, http://www.ce.berkeley.edu/projects/feap/.

\bibitem{art:AinCoy2003}
{\sc M.~Ainsworth and J.~Coyle}, {\em Hierarchic finite element bases on
  unstructured tetrahedral meshes}, Internat. J. Numer. Methods Engrg., 58
  (2003), pp.~2103--2130.

\bibitem{AlnLog08}
{\sc M.~Aln\ae{}s, A.~Logg, K.-A. Mardal, O.~Skavhaug, and H.~P. Langtangen},
  {\em {UFC} Specification and User Manual} 1.1,
  http://www.fenics.org/ufc/ (2008).

\bibitem{art:ArnBofFal2005}
{\sc D.~N. Arnold, D.~Boffi, and R.~S. Falk}, {\em Quadrilateral {$H({\rm
  div})$} finite elements}, SIAM J. Numer. Anal., 42 (2005), pp.~2429--2451.

\bibitem{ArnFalWin06}
{\sc D.~N. Arnold, R.~S. Falk, and R.~Winther}, {\em Finite element exterior
  calculus, homological techniques, and applications}, Acta Numer., 15 (2006),
  pp.~1--155.

\bibitem{ArnFalWin07}
{\sc D.~N. Arnold, R.~S. Falk, and R.~Winther}, {\em Mixed finite element
  methods for elasticity with weakly imposed symmetry}, Math.
  Comp., 76 (2007), pp.~1699--1723.

\bibitem{ArnWin02}
{\sc D.~N. Arnold and R.~Winther}, {\em Mixed finite elements for elasticity},
  Numer. Math., 92 (2002), pp.~401--419.

\bibitem{www:deal.II}
{\sc W.~Bangerth, R.~Hartmann, and G.~Kanschat}, {\em  deal.{\rm{II}}
  {D}ifferential {E}quations {A}nalysis {L}ibrary}, http://www.dealii.org/ (2006).

\bibitem{BreDou87}
{\sc F.~Brezzi, J.~Douglas, Jr., M.~Fortin, and L.~D. Marini}, {\em Efficient
  rectangular mixed finite elements in two and three space variables}, RAIRO
  Mod\'el. Math. Anal. Num\'er., 21 (1987), pp.~581--604.

\bibitem{BreDou85}
{\sc F.~Brezzi, J.~Douglas, Jr., and L.~D. Marini}, {\em Two families of mixed
  finite elements for second order elliptic problems}, Numer. Math., 47 (1985),
  pp.~217--235.

\bibitem{BreFor91}
{\sc F.~Brezzi and M.~Fortin}, {\em Mixed and Hybrid Finite Element Methods},
  Springer Ser. Comput. Math. 15, Springer-Verlag, New
  York, 1991.

\bibitem{CasRie05}
{\sc P.~Castillo, R.~Rieben, and D.~White}, {\em FEMSTER: An object-oriented
  class library of high-order discrete differential forms}, ACM Trans. Math.
  Software, 31 (2005), pp.~425--457.

\bibitem{Cia76}
{\sc P.~G. Ciarlet}, {\em Numerical Analysis of the Finite Element
Method}, S\'{e}min. Math. Sup\'{e}r. 59, Les Presses de
l'Universit\'{e} de Montr\'{e}al, Montreal, 1976.

\bibitem{DomSay08}
{\sc V.~Dom\'{i}nguez and F.-J. Sayas}, {\em Algorithm $884${\rm:} A simple {Matlab}
  implementation of the {Argyris} element}, {ACM} Trans. Math.
  Software, 35 (2009), 11 pp.

\bibitem{logg:www:03}
{\sc J.~Hoffman, J.~Jansson, C.~Johnson, M.~G. Knepley, R.~C. Kirby, A.~Logg,
  L.~R. Scott, and G.~N. Wells}, {\em {FE}ni{CS}}, http://www.fenics.org/ (2006).

\bibitem{Kir04}
{\sc R.~C. Kirby}, {\em Algorithm $839${\rm:} {FIAT}, a new paradigm for computing
  finite element basis functions}, {ACM} Trans. Math. Software,
  30 (2004), pp.~502--516.

\bibitem{www:FIAT}
{\sc R.~C. Kirby}, {\em {FIAT}},
http://www.fenics.org/fiat/ (2006).

\bibitem{Kir06}
{\sc R.~C. Kirby}, {\em Optimizing {FIAT} with {Level $3$ BLAS}}, {ACM} Trans. Math. Software, 32 (2006), pp.~223--235.

\bibitem{KirKne2005}
{\sc R.~C. Kirby, M. Knepley, A.~Logg, and L.~R. Scott}, {\em Optimizing the
  evaluation of finite element matrices}, SIAM J. Sci. Comput., 27 (2005),
  pp.~741--758.

\bibitem{KirLog2006}
{\sc R.~C. Kirby and A.~Logg}, {\em A compiler for variational forms}, {ACM}
  Trans. Math. Software, 32 (2006), pp.~417--444.

\bibitem{KirLog2007a}
{\sc R.~C. Kirby and A.~Logg}, {\em Efficient compilation of a class
of variational forms}, {ACM} Trans. Math. Software, 33 (2007), 20 pp.

\bibitem{Kirby:2008:BDC}
{\sc R.~C. Kirby and A.~Logg}, {\em Benchmarking
  domain-specific compiler optimizations for variational forms}, {ACM}
  Trans. Math. Software, 35 (2008), 18 pp.

\bibitem{KirLogEtAl2006}
{\sc R.~C. Kirby, A.~Logg, L.~R. Scott, and A.~R. Terrel}, {\em Topological
  optimization of the evaluation of finite element matrices}, {SIAM} J. Sci.
  Comput., 28 (2006), pp.~224--240.

\bibitem{KirSco07}
{\sc R.~C. Kirby and L.~R. Scott}, {\em Geometric optimization of the
  evaluation of finite element matrices}, SIAM J. Sci. Comput., 29 (2007),
  pp.~827--841.

\bibitem{logg:www:04}
{\sc A.~Logg}, {\em {FFC}}, http://www.fenics.org/ffc/ (2006).

\bibitem{Log2007a}
{\sc A.~Logg}, {\em Automating the finite element method}, Arch. Comput. Methods Eng., 14 (2007), pp.~93--138.

\bibitem{MarTai02}
{\sc K.~A. Mardal, X.-C. Tai, and R.~Winther}, {\em A robust finite element method for Darcy--Stokes flow}, SIAM J. Numer. Anal., 40 (2002), pp.~1605--1631.

\bibitem{Mon03}
{\sc P.~Monk}, {\em Finite Element Methods for Maxwell's Equations}, Oxford
 University Press, New York, 2003.

\bibitem{Ned80}
{\sc J.-C. N{\'e}d{\'e}lec}, {\em Mixed finite elements in {${\bf R}\sp{3}$}}, Numer. Math., 35 (1980), pp.~315--341.

\bibitem{Ned86}
{\sc J.-C. N{\'e}d{\'e}lec}, {\em New mixed finite elements in {${\bf R}\sp{3}$}}, Numer. Math., 50 (1986), pp.~57--81.

\bibitem{OelLog2008a}
{\sc K.~B.~\O lgaard, A.~Logg, and G.~N. Wells}, {\em Automated
code generation
  for discontinuous Galerkin methods}, SIAM J. Sci. Comput., 31 (2008),
  pp.~849--864.

\bibitem{OelWel09}
{\sc K.~B.~\O lgaard and G.~N. Wells}, {\em Optimisations for quadrature representations of
finite element
  tensors through automated code generation}, ACM Trans. Math. Software, to appear.

\bibitem{www:FreeFEM}
{\sc O.~Pironneau, F.~Hecht, A.~L. Hyaric, and K.~Ohtsuka}, {\em
Free{FEM}}, http://www.\allowbreak freefem.\allowbreak org/ (2006).

\bibitem{RavTho77b}
{\sc P.-A. Raviart and J.~M. Thomas}, {\em Primal hybrid finite element methods
  for {$2$}nd order elliptic equations}, Math. Comp., 31 (1977), pp.~391--413.

\bibitem{Sch03}
{\sc A.~Schneebeli}, {\em An ${H}(\mathrm{curl}; \omega)$ Conforming {FEM}:
  {N}\'ed\'elec's Element of the First Type}, Technical report, 2003; also available
  online from http://www.dealii.org/\allowbreak developer/\allowbreak reports/ \allowbreak nedelec/nedelec.pdf.

\bibitem{www:ngsolve}
{\sc J.~Sch\"oberl}, {\em {NGS}olve},
http://www.hpfem.jku.at/ngsolve/index.html/ (2008).

\end{thebibliography}
\end{document}